\documentclass[12pt,leqno]{amsart}
\theoremstyle{plain}

\textheight9in \def\DATE{\today}
\newtheorem{theorem}{Theorem}
\newtheorem{definition}[theorem]{Definition}

\newtheorem{example}[theorem]{Example}
\newtheorem{lemma}[theorem]{Lemma}
\newtheorem{claim}[theorem]{Claim}

\newtheorem{proposition}[theorem]{Proposition}
\newtheorem{remark}[theorem]{Remark}
\newtheorem{remarks}[theorem]{Remarks}

\def\ext{\mbox{\large$\land$}}
\def\mod{{\it mod}}
\def\Im{{\it Im}}
\def\Ker{{\it Ker}}
\def\span{{\it Span}}

\def\pt{{\bf pt}}
\long\def\comment#1\endcomment{}

\def\endex{\hfill\rule{10pt}{3pt}}

\def\bbN{{\mathbb N}}
\def\bbQ{{\mathbb Q}}
\def\bbK{{\mathbb K}}
\def\bbR{{\mathbb R}}
\def\bbC{{\mathbb C}}
\def\bbP{{\mathbb P}}

\def\kk{\kappa }
\def\otexp#1#2{{#1}^{\otimes #2}}
\def\ot{\otimes}
\def\rada#1#2{#1,\ldots,#2}
\def\Rada#1#2#3{#1_{#2},\dots,#1_{#3}}
\def\DD#1#2{{\Delta_{#1,#2}}}

\def\fatDelta{{\hbox{\hskip 1pt$\Delta \hskip -10.5pt%
                \Delta \hskip -10.5pt \Delta\hskip 1pt$}}}
\def\oH{{%
   \hskip 3.235pt\raisebox{8.7pt}{\scriptsize\rm o \hskip -4pt} \hskip
-7pt H}}
\def\lloH{{%
   \hskip 3.235pt\raisebox{6.7pt}{\scriptsize\rm o \hskip -4pt} \hskip
-7pt H}}
\def\cases#1#2#3#4{
                  \left\{
                         \begin{array}{ll}
                           #1,\ &\mbox{#2}
                           \\
                           #3,\ &\mbox{#4}
                          \end{array}
                   \right.
}
\def\oHsuper#1{{\oH^{\hskip .5pt\raisebox{-4.4pt}{\scriptsize $#1$}}}}
\def\Kriz{K\v r\' \i \v z}\def\kriz{\Kriz}
\def\Ker{{\it Ker}}
\def\ext{\mbox{\large$\land$}}

\def\oD{{%
   \hskip 3.235pt\raisebox{8.7pt}{\scriptsize\rm o \hskip -3pt} \hskip
-7pt \Delta}}
\def\ssoD{{%
   \hskip 3.235pt\raisebox{6pt}{\tiny\rm o \hskip -3pt} \hskip
-5pt \Delta}}
\def\oX{{%
   \hskip 3.235pt\raisebox{8.7pt}{\scriptsize\rm o \hskip -4pt} \hskip
-7pt X}}
\def\lloX{{%
   \hskip 3.235pt\raisebox{6.7pt}{\scriptsize\rm o \hskip -4pt} \hskip
-7pt X}}
\def\ssoX{{%
   \hskip 2.235pt\raisebox{6pt}{\tiny\rm o \hskip -4pt} \hskip
-5pt X}}
\def\od{{%
   \hskip 3.235pt\raisebox{8.7pt}{\scriptsize\rm o \hskip -4pt} \hskip
-7pt d}}
\def\upperod{{%
   \hskip 3.235pt\raisebox{6.7pt}{\tiny\rm o \hskip -4pt} \hskip
-4pt d}}
\def\oK{{%
   \hskip 3.235pt\raisebox{8.7pt}{\scriptsize\rm o \hskip -4pt} \hskip
-7pt K}}
\def\oY{{%
   \hskip 3.235pt\raisebox{8.7pt}{\scriptsize\rm o \hskip -3.3pt} \hskip
-7pt Y}}
\def\Ann{{\rm Ann}}
\def\AA{{\mathcal A}}
\def\BB{{\mathcal B}}
\def\CC{{\mathcal C}}
\def\FF{{\mathcal F}}
\def\EE{{\mathcal E}}
\def\HH{{\mathcal H}}
\def\KK{{\mathcal K}}
\def\MM{{\mathcal M}}
\def\NN{{\mathcal N}}
\def\TT{{\mathcal T}}
\def\cinf{C^{\infty }}

\def\wtY{\widetilde{Y}}
\def\wtZ{\widetilde{Z}}
\def\bE{{\bf E}}
\def\ss{\! \mbox{ }_W\bE _1}
\def\spos{\! \mbox{ }_W\bE _0}
\def\s2{\! \mbox{ }_W\bE _2}

\def\sir{\! \mbox{ }_W\bE _{r}}
\def\dR{\Omega^*_{dR}}
\def\harp{
     \unitlength=.4cm
          \begin{picture}(1.3,1)(-0.1,0.15)
            \thicklines
            \put(0.00,0.3){\line(1,0){1}}
            \put(0.00,0.7){\line(1,0){1}}
            \put(0.3,0.00){\line(0,1){1}}
            \put(0.7,0.00){\line(0,1){1}}
          \end{picture}}

\title[Multiplicative models for configuration spaces]%
      {Multiplicative models for configuration spaces of algebraic varieties}

\author[B. Berceanu]{Barbu Berceanu}
\author[M. Markl]{Martin Markl}
\author[\c S. Papadima]{\c Stefan Papadima}
\thanks{The second author supported by an EURROMMAT
Invited Lecturer Fellowship (Contract ICA1--CT--2000--70022)
and by the grant GA AV {\v CR}
201/02/1390. The first and third author
partially supported by CNCSIS grant 693/2002 of the Romanian Ministry
of Education and Research.}

\begin{document}

\bibliographystyle{plain}
\baselineskip16pt plus 1pt minus .5pt

\begin{abstract}
\baselineskip13.5pt
\parskip3pt
W.~Fulton--R.~MacPherson~\cite{FMP} found a Sullivan dg-algebra model
for the space of $n$-configurations of a smooth compact complex algebraic
variety $X$. I.~\kriz~\cite{kriz} gave a simpler model, $E_n(H)$,
depending only on the cohomology ring, $H:=H^*X$.

We construct an even simpler and smaller model, $J_n(H)$. We then
define another new dg-algebra, $E_n(\lloH)$, and use $J_n(H)$ to prove
that $E_n(\lloH)$ is a model of the space of $n$-configurations of the
{\em non-compact} punctured manifold $\lloX$, when $X$ is
$1$-connected. Following an idea of V.G.~Drinfel'd \cite{D}, we put a
simplicial bigraded differential algebra structure on~$\{E_n(\lloH)
\}_{n \geq 0}$.
\end{abstract}

\maketitle

\section{\bf Introduction and the main results}
\label{intro}

Let $X$ be a connected space.
The topology of {\em ordered configuration spaces}
\[
F(X,n) := \{(\Rada x1n) \in X^{\times n};\ x_i \not= x_j
\mbox { if } i \not= j\}
\]
of $n$ distinct labeled points in $X$
has attracted a considerable attention, over the years.

The cohomology rings $ H^*F(\bbR ^2,n) $ have been described by
Arnold~\cite{Ar}. In his 1972 thesis, F.~Cohen extended Arnold's
computations to all Euclidean spaces; see~\cite{C}.
For $ X $ an $ l $-dimensional real oriented manifold, the Leray
spectral sequence of the inclusion $ F(X,n)\hookrightarrow X^{\times n} $
has been described by Cohen-Taylor~\cite{CT} and further analyzed by
Totaro~\cite{totaro}. With field coefficients $ \bbK $, the above
Cohen-Taylor spectral sequence converges multiplicatively to
$ H^*(F(X,n); \bbK) $; it has the property that $ \bE _2=\bE _l $, and
the differential graded algebra $ (\bE _l,d_l) $ depends only on $ n $
and the cohomology algebra $ H^*(X; \bbK) $. See~\cite{CT},~\cite{totaro}.
If $X$ is a smooth projective complex $ m $-variety, Totaro~\cite{totaro}
showed, over $ \bbQ $, that $ {\bf E}_{2m+1}={\bf E}_{\infty } $, and
$ H^*(F(X,n); \bbQ)=\bE _{2m+1} $, as graded algebras.

For a {\em compact} oriented real $ l $-manifold $ X $, it is convenient
to view $ F(X,n) $ as $ X^{\times n}\setminus D^n X $, where $ D^nX $
denotes the fat diagonal. Using Lefschetz duality, one may thus replace
the Betti numbers of $ F(X,n) $ by those of the pair
$ (X^{\times n}, D^n X) $, modulo suitable regrading. In this way, Brown
and White~\cite{BW} were able to compute the Betti numbers of $ F(X,n) $
in terms of the cohomology algebra of $ X $, for $ n\leq 3 $. For
arbitrary $ n $, Bendersky and Gitler have constructed in~\cite{BG}
another spectral sequence, converging {\em additively} to
$ H^*(F(X,n); \bbK) $, regraded via Lefschetz duality. They have also
proved that $ \bE _2=\bE _\infty $ in their spectral sequence over $ \bbQ $,
when $ X $ is rationally formal, in the sense of Sullivan~\cite{S}
(for instance, when $ X $ is compact K\"ahler manifold,
see~\cite{DGMS}).

As it turns out, there is an additive isomorphism between the
Bendersky-Gitler $ \bE _2 $-term and the Cohen-Taylor $ \bE _{l+1} $-term,
after regrading (this was proved independently in~\cite[Theorem~29]{MP},
and in~\cite[Theorem~1]{FT}). Therefore, the additive part of Totaro's
collapsing result actually holds for $ \bbQ $-formal closed oriented
manifolds. Nevertheless, the above spectral sequences for $ F(X,4) $
do not collapse in general, as follows from the example given by F\' elix
and Thomas~\cite{FT}, where $X$ is the sphere tangent bundle of
$ S^2\times S^2 $.

Our aim in this paper is to go beyond Betti numbers and cohomology algebras.
We will describe two new {\em differential graded algebra (DGA) models}, for
$ \bbK $-homotopy types in characteristic zero (in the sense of
Sullivan~\cite{S}) of configuration spaces, $ F(X,n) $ and $ F(\oX,n) $.
Here $ X $ is a smooth compact complex algebraic $m$-variety, and
$ \oX:=X\setminus \pt $ denotes the punctured manifold.

Let $Y$ be a smooth complex algebraic variety. It is known that $Y$
has a convenient {\rm compactification}. That is, $ Y=\wtY\setminus D
$, where $ \wtY $ is smooth compact, and $ D\subset \wtY $ is a
divisor with normal crossings. A basic result of
Morgan~\cite[Theorem~9.6 and Corollary~9.7]{M} says then that the $
\bbC $-homotopy type of $Y$ is naturally determined by the $ \bbC
$-cohomology algebras of the various intersections of components of
$D$, together with the restriction and Gysin maps between them. The
same holds true (non-naturally) over $\bbQ$\hskip .5mm;
see~\cite[Theorem~10.1]{M}.

Consider now $ Y=F(X,n) $, where $X$ is a smooth compact complex algebraic
variety. In a seminal paper, Fulton and MacPherson~\cite{FMP}
constructed a particularly nice compactification $ \wtY $ of this space,
in the above sense. Applying Morgan's
theory, they succeeded to describe a DGA model for the $ \bbQ $-homotopy
type of $ F(X,n) $, depending on $n$, $ H^*(X; \bbQ )$, and the Chern
classes of $X$. This model was algebraically simplified by \kriz~\cite{kriz},
whose $ \bbQ $-model for $F(X,n)$ depends only on $n$ and the algebra
 $ H^*(X; \bbQ )$. To describe our models, we begin by introducing a
construction which abstracts the key features of the \Kriz \/ model.

Let $A$ be a unital graded commutative algebra over a field $ \bbK $,
 $m \geq 1$ a fixed number
and $\nabla \in A \ot A$ a degree $2m$ class which is graded
symmetric:
\begin{equation}
\label{eli1}
T(\nabla) = \nabla,
\end{equation}
where $T(a \ot b) := (-1)^{\deg(a)\deg(b)}(b \ot a)$ is the graded flip,
and which is `diagonal' in the sense that
\begin{equation}
\label{eli2}
(a \ot 1)\nabla = (1 \ot a)\nabla,
\end{equation}
for any $a \in A$. For example, if $ A=H^*(X;\bbK ) $, where $X$ is an
oriented $ 2m $-dimensional real manifold, then the diagonal class
$ \nabla \in H^{2m}(X\times X;\bbK ) $ will satisfy (1) and (2) above.
Given a subset $\{ i_1< \ldots <i_k\} \subset \{\rada 1n\}$,
denote by
\[
\iota_{\Rada i1k} : \otexp Ak \hookrightarrow \otexp An
\]
the obvious inclusion that puts the $s$-th factor of $\otexp Ak$ into
the $i_{s}$-th slot of $\otexp An$. More generally, for an arbitrary, not
necessarily linearly ordered subset
$\{j_1,\ldots,j_k\} \subset \{\rada 1n\}$, there is a unique
permutation $\sigma \in \Sigma_k$ such that $i_s := j_{\sigma^{-1}(s)}$,
$1 \leq s \leq k$, is linearly ordered. We will put
\[
\iota_{\Rada j1k} : = \iota_{\Rada i1k} \circ  T_\sigma,
\]
where $T_\sigma : \otexp Ak \to \otexp Ak$ is
the canonical automorphism induced by $\sigma$.
Later in the paper we will also need
\begin{equation}
\label{eli7}
\iota_{\varphi} := \iota_{\varphi (1),\ldots, \varphi (k)}:
\otexp Ak \hookrightarrow \otexp An,
\end{equation}
where  $\varphi : \{ \rada 1k \}\to \{ \rada 1n \}$ is an
injective map.
Set $\nabla_{ij} := \iota_{ij}(\nabla) \in \otexp An$, for $ 1
\leq i \neq j \leq n$.

\begin{definition}
\label{eli}
Let $E_n(A,\nabla)$ be the free graded
commutative $\otexp An$-algebra
$\otexp An[G_{ij}]$ with degree $2m-1$ exterior generators $G_{ij}$, $n \geq
i > j \geq 1$, modulo the following relations:
\begin{eqnarray}
\label{one}
G_{ij} G_{ik}  &=& G_{jk}(G_{ik} - G_{ij}),\ \mbox { for }
n \geq i > j > k \geq 1, \mbox { and}
\\
\label{2}
 \iota_i(x)G_{ij} &=& \iota_j(x)G_{ij},\mbox { for $x \in A$, $n \geq
i > j \geq 1$,}
\end{eqnarray}
with the differential $d$ given by $d(\iota_i(x)) := 0$ for  $n \geq i
\geq 1$ and $d G_{ij} := \nabla_{ij}$ for $n \geq i > j \geq 1$.
\end{definition}

We leave to the reader to verify that symmetry~(\ref{eli1}) guarantees the
compatibility of $d$ with relations~(\ref{one}), while~(\ref{eli2}) is
necessary for $d$ to be compatible with~(\ref{2}).

{}From now on, $H$ will denote
an arbitrary $2m$-dimensional Poincar{\'e} duality algebra over a field
$\bbK$, with a distinguished orientation class
$\omega \in H^{2m}\setminus \{ 0\}$. Let $\Delta \in H^{\ot 2}$  be
the class of the diagonal,
\[
\Delta := \sum_{\alpha} (-1)^{\deg (h_{\alpha})} h_{\alpha}
  \otimes h^*_{\alpha},
\]
where $\{ h_{\alpha} \}$ is a
homogeneous $\bbK$-basis of $H^*$, comprising $1$, and
$\{ h^*_{\alpha} \}$ is the Poincar\' e dual basis with respect to
$\omega$. The couple $(H,\Delta)$ then clearly fulfills
both~(\ref{eli1}) and~(\ref{eli2}).
The {\em \Kriz-model\/} $E_n(H)$ is defined by
$E_n(H) : = E_n(H,\Delta)$. When $ H=H^*(X;\bbK ) $, where $X$ is a closed
oriented real $ 2m $-manifold, note also that $ (E_n(H),d) $ coincides
with the DGA $ (\bE _{2m},d_{2m}) $ coming from the Cohen-Taylor spectral
sequence for $ F(X,n) $; see~\cite{totaro}.

\begin{theorem}
[Fulton-MacPherson~\cite{FMP}, K\v r\' \i \v z~\cite{kriz}]
\label{thm:basic}
Let $X$ be a smooth compact complex algebraic variety and $
H=H^*(X;\bbQ )$.  Then the {\em DGA} $(E_n(H),d)$ is a rational model,
in the sense of Sullivan, of the configuration space $F(X,n)$.
\end{theorem}

Let $(\fatDelta^n) := (\DD 21,\DD
31,\ldots,\DD n1)$ denote the ideal generated in $\otexp Hn$ by
the diagonals $\DD 21,\DD 31,\ldots,\DD n1$. Let $\pi :
\otexp Hn \to
{\otexp  Hn}/{(\fatDelta^n)}$
be the projection and $\jmath_i : H \to {\otexp  Hn}/{(\fatDelta^n)}$
the composition $\pi \circ \iota_i$, $n \geq i
\geq 1$. Our first DGA model, $J_n(H)$, is basically
the \kriz-model $E_n(H)$
with $\otexp Hn$ replaced by $\otexp Hn/(\fatDelta^n)$ and $n-1$
generators, $G_{n1},G_{n-1,1},\ldots,G_{21}$, missing:

\begin{definition}
\label{def:j}
The model $J_n(H)$ is defined to be the free graded commutative
$\otexp Hn/(\fatDelta^n)$-algebra
$\otexp Hn/(\fatDelta^n)[G_{ij}]$
with degree $2m-1$ exterior generators $G_{ij}$, $n \geq i
> j \geq 2$, modulo the relations:
\begin{eqnarray*}
G_{ij} G_{ik}  &=& G_{jk}(G_{ik} - G_{ij}),\
n \geq i > j > k \geq 2, \mbox { and}
\\
 \jmath_i(x) G_{ij} &=&  \jmath_j(x) G_{ij} ,\mbox { for $x \in H$,
$n \geq i > j  \geq 2$.}
\end{eqnarray*}
The differential is determined by $d(\jmath_i(x)): = 0$ for $n
\geq i \geq 1$ and $d G_{ij} := \pi(\Delta_{ij})$ for $n \geq i > j
\geq 2$.
\end{definition}

Observe that $J_n(H)$ is not of the form $E_n(A,\nabla)$, but it is
very close to it.
Though the `coefficients' $\otexp Hn/(\fatDelta^n)$ of $J_n(H)$ seem
more complicated than the coefficients $\otexp Hn$ of the \kriz-model,
we will see, in Proposition~\ref{oH}, that the structure of
$\otexp Hn/(\fatDelta^n)$ is actually very simple.
The projection $\pi$ clearly generates an epimorphism
$\Psi : E_n(H) \to J_n(H)$ of DG-algebras
with $\Psi(G_{ij}) :=  G_{ij}$ if $n \geq i>j  \geq 2$ and
$\Psi(G_{j1}) := 0$, $n \geq j \geq 2$.
Our first main result in this paper, whose proof we postpone to
Section~\ref{proofs}, reads:

\begin{theorem}
\label{main}
The natural map $\Psi : E_n(H) \to J_n(H)$ induces a cohomology
isomorphism, for an arbitrary even-dimensional Poincar\' e duality algebra
$H$. Therefore, $J_n(H)$ is also a {\em DGA}
model, in the sense of Sullivan, of the configuration space $F(X,n)$, for
$X$ and $H$ as in Theorem~\ref{thm:basic}.
\end{theorem}

Let $H$ be an arbitrary even-dimensional Poincar\' e duality algebra,
 as before. We are going to consider another associated DGA, to be denoted
by $E_n(\oH)$. To begin with, let
$\oH$ be the quotient algebra, $\oH :=H/\bbK \cdot \omega$,
that is,
\[
\oHsuper i
= \cases{H^i}{for $0 \leq i < 2m$\, , and}{0}{for $i\geq 2m$ \, ,}
\]
with multiplication induced from $H$. Notice that,
when $ H=H^*(X;\bbK )$, with $X$
a closed oriented real $ 2m $-manifold, $\oH$ is nothing
else but the cohomology algebra of the non-compact
punctured manifold $\oX = X \setminus \pt $.
Denote by $\oD$ the image of $\Delta$ in $\oH^{\ot 2}$.
Plainly, conditions (1) and (2) are satisfied by $ (\oH ,\oD ) $.

\begin{definition}
\label{def:circ}
The {\em punctured \kriz -model} is
the differential graded commutative algebra $E_n(\oH) := E_n(\oH,\oD)$,
with differential denoted by $ \od $.
\end{definition}

There is a very intimate relation between the $J$-model from
Definition~\ref{def:j} and the punctured \kriz-model, which can be
described as follows.  Observe first that $H$ is naturally augmented,
via an augmentation $\epsilon : H \to \bbK$, which makes $\bbK$ a
right $H$-module.  Observe also that the map $\iota_1 : H \to \otexp
Hn$ induces a left differential $ (H,0) $-module structure on
$J_n(H)$.  The following useful fact will be proved in Section 3:

\begin{proposition}
\label{eli3}
For any $n \geq 1$, there is an isomorphism of differential graded
commutative algebras
\[
(E_{n-1}(\oH),\od ) \cong \bbK \ot_{(H,0)}(J_n(H),d)\, .
\]
\end{proposition}

In Section~\ref{sec:punct},
we will derive  from Theorem~\ref{main} and Proposition~\ref{eli3}
the second main result of our paper, which may be viewed as
an analog of the fundamental Theorem~\ref{thm:basic} in a
{\em non-compact } situation.

\begin{theorem}
\label{thm:punct}
Let $X$ be a $1$-connected  smooth compact complex algebraic variety.
Set $ H=H^*(X;\bbC ) $. Then
the differential graded algebra $(E_n(\oH),\od) $ is a $\bbC$-model,
in the sense of Sullivan, of
the configuration space $F(X\setminus \pt,n)$.
\end{theorem}

Our proof entails a careful analysis of the natural fibration
$$F(\oX ,n-1)\hookrightarrow  F(X,n)\stackrel p\longrightarrow  X \ni \pt .$$
Here $p$ is the projection onto the first coordinate, which is induced by an
algebraic map defined on the Fulton-MacPherson compactification of
$ F(X,n) $, see~\cite{FMP}. A key step involves naturality properties
from Morgan's theory~\cite{M}, which explains our need to use $ \bbC $
instead of $ \bbQ $ coefficients.

The strength
of  Theorem~\ref{thm:punct} is illustrated by Example~\ref{ex:OS}: in the
simplest case, corresponding to $ X=\bbC \bbP ^1 $, we
very easily recover the formality of the classifying space
$F(\bbR^2,n)$ of the pure braid group, as well as Arnold's description
of its cohomology.

The rich geometry of configuration spaces of manifolds has a natural
algebraic analogue, at the level of $E$-models. We will illustrate
this principle in \S~\ref{sec:struct}. In Proposition~\ref{prop:drin},
we endow the collection
$\{ E_n(\oH) \}_{n \geq 0}$ with a simplicial structure, in the
category of {\em bigraded differential algebras\/} (DBGA's); this is based on
Drinfel'd's~\cite{D} cosimplicial group structure on Artin pure braid
groups, and on the doubling operations on chord diagrams from the
theory of Vassiliev invariants of links~\cite{BN}. We also define a coaction
map, relating the cohomology of little cubes to
   $\{ E_n(\oH) \}_{n \geq 0}$, in
Proposition~\ref{prop:coact}. It  should be pointed out that both
above-mentioned extra structures on DGA models exist {\em only\/} for
punctured models; see Remarks~\ref{rks:pess1} and~\ref{rk:pess2}.

\section{\bf Hints for applications}
\label{sec:hints}

The algebra $E_n(A,\nabla)$ introduced in Definition~\ref{eli}
is not free as a graded commutative
algebra, but it can be presented as a direct sum of free $\otexp
A{n-k}$-modules, $0 \leq k \leq n-1$.
To formulate this statement more
precisely, we need the following notation.

For a sequence $1 \leq i_1 < \cdots < i_k \leq n$,
let $1 \leq h_1 < \cdots < h_{n-k} \leq n$ be its `complement', that
is, a sequence such that
$\{\Rada i1k,\Rada h1{n-k}\} = \{\rada 1n\}$. Set
\[
\iota^c_{\Rada i1k} : = \imath_{\Rada h1{n-k}} : \otexp A{n-k} \to
\otexp An.
\]
The following statement is Proposition~2.1 of~\cite{Bez}.

\begin{proposition}
[Bezrukavnikov~\cite{Bez}]
\label{bezru}
The linear map
\begin{equation}
\label{22}
\Xi:
\bigoplus_{0 \leq k \leq n-1}\bigoplus_{I_k}
\otexp A{n-k}\cdot G_{i_1j_1} G_{i_2j_2} \cdots G_{i_kj_k}
\longrightarrow E_n(A),
\end{equation}
where $I_k := \{2 \leq i_1 < i_2 < \cdots < i_k \leq n,\
i_1 > j_1 \geq 1, \ldots ,i_k > j_k \geq 1\}$,
given by
\[
\Xi(h \cdot  G_{i_1j_1} G_{i_2j_2} \cdots
G_{i_kj_k}) := \iota^c_{\Rada i1k}(h) \cdot
G_{i_1j_1} G_{i_2j_2} \cdots
G_{i_kj_k},
\]
is an isomorphism of graded vector spaces.
\end{proposition}

Let now $H$ be an even-dimensional Poincar\' e duality algebra over a
field $ \bbK $. The direct sum~(\ref{22}) with the
induced differential can be therefore understood as an alternative
description of $E_n(H)$ or $E_n(\oH)$ (depending on whether
$(A,\nabla)$  is $(H,\Delta)$ or $(\oH,\oD)$),
accessible by methods of linear algebra.

\begin{example}
\label{full}
{\rm\
For $n=3$, the direct sum decomposition~(\ref{22}) of $E_3(H)$ equals:

\begin{center}
{
\unitlength=1pt
\begin{picture}(180.00,123.00)(0.00,0.00)
\put(140.00,5.00){\makebox(0.00,0.00){\scriptsize $\Delta_{21}$}}
\put(140.00,55.00){\makebox(0.00,0.00){\scriptsize $-\Delta_{32}$}}
\put(140.00,82.00){\makebox(0.00,0.00){\scriptsize $-\Delta_{31}$}}
\put(140.00,115.00){\makebox(0.00,0.00){\scriptsize $\Delta_{21}$}}
\put(30.00,20.00){\makebox(0.00,0.00){\scriptsize $\Delta_{32}$}}
\put(40.00,70.00){\makebox(0.00,0.00){\scriptsize $\Delta_{21}$}}
\put(30.00,102.00){\makebox(0.00,0.00){\scriptsize $\Delta_{31}$}}
\put(200.00,30.00){\makebox(0.00,0.00){$H \cdot G_{21} G_{32}$}}
\put(200.00,90.00){\makebox(0.00,0.00){$H \cdot G_{21} G_{31}$}}
\put(160.00,20.00){\vector(-3,-1){52}}
\put(160.00,40.00){\vector(-3,1){52}}
\put(160.00,80.00){\vector(-3,-1){52}}
\put(160.00,100.00){\vector(-3,1){52}}
\put(80.00,0.00){\makebox(0.00,0.00){$\otexp H2 \cdot G_{32}$}}
\put(80.00,60.00){\makebox(0.00,0.00){$\otexp H2 \cdot G_{21}$}}
\put(80.00,120.00){\makebox(0.00,0.00){$\otexp H2 \cdot G_{31}$}}
\put(52,8){\vector(-1,1){42}}
\put(52,60.00){\vector(-1,0){35}}
\put(52,112){\vector(-1,-1){42}}
\put(0.00,60.00){\makebox(0.00,0.00){$\otexp H3$}}
\end{picture}}
\end{center}

\noindent
where the maps are multiplications by indicated elements (modulo Koszul signs).
}
\endex
\end{example}

One can easily adapt Bezrukavnikov's proof to obtain
 the following analog of decomposition~(\ref{22})
also for $J_n(H)$.

\begin{proposition}
\label{ifgod}
The linear map
\begin{equation}
\label{alt}
\Upsilon:
\bigoplus_{0 \leq k \leq n-2}
\bigoplus_{J_k}
\otexp H{n-k}/{(\fatDelta^{n-k})}
\cdot G_{i_1j_1} G_{i_2j_2}\cdots G_{i_kj_k}
\longrightarrow J_n(H),
\end{equation}
where $ J_k := \{3 \leq i_1 < i_2 <
\cdots < i_k \leq n,\ i_1 > j_1 \geq 2, \ldots ,i_k > j_k \geq 2\}$,
given by
\[
\Upsilon(h \cdot  G_{i_1j_1} G_{i_2j_2} \cdots
G_{i_kj_k}) := \jmath^c_{\Rada i1k}(h) \cdot
G_{i_1j_1} G_{i_2j_2} \cdots G_{i_kj_k}
\]
is an isomorphism of graded vector spaces. In the above display,
$\jmath^c_{\Rada i1k}: \otexp H{n-k}/(\fatDelta^{n-k}) \to \otexp
Hn/(\fatDelta^n)$ is the map induced by $\iota^c_{\Rada i1k}$.
\end{proposition}

The left hand side of~(\ref{alt}) provides an alternative description
of
$J_n(H)$ as a direct sum of free $\otexp
H{n-k}/{(\fatDelta^{n-k})}$-modules.
{}From this perspective, the following proposition (to be proved in
Section~\ref{proofs}) is very useful.

\begin{proposition}
\label{oH}
The composition
\[
H \ot \oHsuper {\ot l-1} \hookrightarrow \otexp Hl
\stackrel\pi\rightarrow \otexp Hl/(\fatDelta^l)
\]
induces, for any $l \geq 2$, an isomorphism of graded left
$H$-modules, where $H$ acts on the first position on both
$H \ot \oHsuper {\ot l-1}$ and $\otexp Hl $.
\end{proposition}

\begin{remark}
\label{bid}
{\rm Note that all three models, $E_n(H)$, $E_n(\oH)$
and $J_n(H)$, are actually {\em bigraded differential algebras\/} (DBGA's). The
second ({\em exterior\/}) degree is given by the number of $G_{ij}$-factors,
and the differentials are homogeneous, of bidegree $(+1,-1)$.
In particular,
\[
H^*F(X ,n)= \bigoplus_{p,q\ge 0} H^{p,q} F(X, n)\, ,
\]
where $p$ denotes the usual degree and $q$ is the exterior degree
(when Theorem~\ref{main} applies), and similarly for $\oX$ (when
Theorem~\ref{thm:punct} applies). We can thus define the bigraded
Poincar\' e polynomial
$$ P_{F(X ,n)}(s,t) :=\sum_{p,q\geq 0}
\mbox {dim}\, H^{p,q}F(X,n)s^q t^p \, . $$
The ordinary Poincar{\'e} polynomial is then the specialization
$P_{F(X ,n)}(t) := P_{F(X ,n)}(1,t)$, and likewise for $ \oX $.
\endex
}
\end{remark}

In Examples 13--17, $X$ is as in Theorem~\ref{main}, and $ H=H^*(X;\bbQ ) $.

\begin{example}
\label{ex:j2}
{\rm\
Let us begin with the simple case of two points (where various methods can
be used). Here $J_2(H)$ is just $\otexp H2/(\Delta)$ with trivial
differential, therefore $F(X,2)$ is a formal space, and
\[
H^*(F(X,2); \bbQ) = \otexp {H}2/(\Delta).
\]

By Proposition~\ref{oH}, the algebras
in the above display are isomorphic, as graded vector spaces, to
$ H\ot \oH $, therefore the Poincar\'e polynomial
$P_{F(X,2)}(s,t)$
of $F(X,2)$ equals $P_X(t)(P_X(t)-t^{2m})$.\endex
}
\end{example}

\begin{example}
\label{ex:j3}
{\rm\
For three points, the model $J_3(H)$ reduces to the 2-term complex
\begin{equation}
\label{nemo}
{\otexp H3}/{(\Delta_{31},\Delta_{21})}
\stackrel d{\longleftarrow} \hskip 1mm
\uparrow^{2m-1} {\otexp H2}/{(\Delta_{21})},
\end{equation}
where $\uparrow^{2m-1}$ denotes the suspension iterated $2m-1$ times
and
the differential $d$ is given by
\[
d(h_1 \ot h_2) = (-1)^{\deg (h_1)+\deg (h_2)}
(h_1 \ot h_2 \ot 1) \cdot \Delta_{32}, \mbox { for
$h_1 \ot h_2 \in \otexp H2$.}
\]
This economical description should be compared with the full
\kriz-model for $F(X,3)$ described in Example~\ref{full}.
Using Proposition~\ref{oH}, we can simplify this complex further to
\[
H \ot \oHsuper{\ot 2}  \stackrel d{\longleftarrow}\hskip 2pt
\uparrow^{2m-1} H \ot \oH ,
\]
where $d$ is now given by
\[
d(h \ot l) = (-1)^{\deg (h)+\deg (l)}[
(h \ot l \ot 1) \Delta_{32} - (h \ot 1 \ot l)\Delta_{21} \!-\!
(-1)^{\deg(h)\deg(l)}(1 \ot l \ot h)\Delta_{31}],
\]
for $ h \ot l \in H \ot \oH $.
As an exercise, check that $ d(h \ot l) $ indeed belongs to $ H \ot
\oHsuper{\ot 2} $. Notice also that the complex~(\ref{nemo}) can be used
to compute the full ring structure of the cohomology of $F(X,3)$, not
only its Betti numbers. \endex
}
\end{example}

\begin{example}
\label{ex:g1}
{\rm
Let $a,b,c$ and $d$ be generators of degree $1$ and $V$ a
two-dimensional graded vector space concentrated in degree $2$
considered as a non-unital algebra with trivial multiplication. Let
$\ext(-)$ denote the free graded commutative associative algebra
functor. Then the rational cohomology algebra of the configuration
space $F(T,3)$ of three points in the two-dimensional
torus $T$ is isomorphic to
\[
\left(
\frac{\ext(a,b,c,d)}{(ab=cd=0,\ ac = bd)} \oplus V
\right)
\ot H^*(T;\bbQ).
\]
The bigraded Poincar{\'e} polynomial $P_{F(T,3)}(s,t)$ is
\[
(1+4t+3t^2+2st^2 )\cdot P_T(t) =1+6t+(12+2s)t^2+
  (10+4s)t^3+ (3+2s)t^4
\]
and the ordinary Poincar\'e polynomial $P_{F(T,3)}(t)$ equals
\hfill\break
\vskip-3pt
\hglue2cm
\hfill $(1+4t+5t^2 )\cdot P_T(t) =1+6t+14t^2+14t^3+ 5t^4.$\hfill \endex\break
}
\end{example}

\begin{example}
\label{ex:gge2}
{\rm
Let $X$ be a Riemann surface of genus $g \geq 2$. Denote by
$\Rada a1g,\Rada b1g$
the standard symplectic basis of $H^1(X;\bbQ)$.
Then the kernel of the map $d$
in~(\ref{nemo}) is the suspension of the ideal generated by
$2g^2 + g$
elements,
\begin{eqnarray*}
&(a_i \ot a_j + a_j \ot a_i), \
(b_i \ot b_j + b_j \ot b_i),\ \mbox {with  $1 \leq i \leq j \leq g$},&
\\
&\mbox {and }\
(a_i \ot b_j + b_j \ot a_i),\ \mbox {with  $1 \leq i,j \leq g$,}&
\end{eqnarray*}
plus one `exceptional' element
\[
1 \ot \omega + 2\sum_1^g (a_i \ot b_j - b_j \ot a_i) - \omega \ot 1.
\]
The
Poincar\'e polynomial of the kernel of $d$ is $t^3(2g^2 +g+1+2gt)$
and the bigraded Poincar\'e polynomial of
$F(X,3)$ equals:
\[
P_{F(X,3)}(s,t) = 1+6gt+12g^2t^2+[8g^3+(2g^2+g+1)s]t^3 +(2g^2+g+2gs)t^4.
\]
Specializing at $s=1$ gives
\[
P_{F(X,3)}(t) = 1+6gt+12g^2t^2+(8g^3+2g^2+g+1)t^3 +(2g^2+3g)t^4,
\]
which is a formula of~\cite{BW}.\endex
}
\end{example}

\begin{example}
\label{ex:g0}
{\rm\
One may use the $J$-model to compute the Poincar\'e polynomial of
$F(S^2,n)$, for arbitrary $n$, in the following simple way. Start
with $ n=3 $, $ J_3(H): H\stackrel{-2\omega}\longleftarrow \hskip 2pt
\uparrow H $ to obtain $ H^*(F(S^2,3))=\ext (x_3) $, where $x_3$
is a degree $3$ generator, with  bigraded Poincar\'e
polynomial $ P_{F(S^2,3)}(s,t)=1+st^3$.
By induction on $ n $, we can find the bigraded Poincar\'e polynomial
$$
P_{F(S^2,n)}(s,t)=(1+st^3)\prod _{k=2}^{n-2}(1+kst) \, ,$$
thus recovering a classical formula, see~\cite[Theorem V.7.1 and Corollary
 V.1.4]{FH}. The induction is based on the following easy argument.
Filter $ J_n(H) $ by $ F_0=J_{n-1}(H) $,  $ F_1=J_n(H) $;
in the corresponding spectral sequence $ \bE = (\bE _n,d_n)$ we obtain
$$ \begin{array}{l}
      (\bE _0,d_0)=(J_{n-1}(H),d)\otimes (\span_\bbQ (1,G_{n2},\dots ,
                         G_{n,n-1}),0) \, , \\
      (\bE _1,d_1)=(H^*(J_{n-1})\otimes (\span_\bbQ (1,G_{n2},\dots ,
                            G_{n,n-1})),0)\, ,
   \end{array} $$
and then $ \bE $ collapses, due to an obvious degree argument. \endex
}
\end{example}

In Examples 18-21, $X$ is as in Theorem~\ref{thm:punct},
and $ H=H^*(X;\bbC ) $.

\begin{example}
\label{ex:OS}
{\rm
Applying Theorem~\ref{thm:punct} to $ X=\bbC \bbP^1 $,
for which $\oH =\bbC \cdot 1$, and $\od =0$, we recover another
classical result: the {\em pure braid space}
$F(\bbR^2 ,n)$ is formal, with cohomology ring described by
the Arnold relations~(\ref{one}), and bigraded Poincar\' e
polynomial given by
\[
P_{F(\bbR^2 ,n)}(s,t) = (1+st)(1+2st)\cdots [1+(n-1)st],
\]
for all $n$; see~\cite{Ar}, and also~\cite{OT}.   \endex
}
\end{example}

\begin{example}
\label{ex:circ12}
{\rm
For $n=1$, Theorem~\ref{thm:punct} reduces to the
statement that $\oX$ is a formal space, a result of
Avramov~\cite{A} (valid for any $1$-connected formal closed
manifold $X$). For $n=2$, it follows from
isomorphisms~(\ref{22}) that the underlying chain complex of our model $E_2(\oH)$ is
\[
\oH^{\ot 2}\stackrel{\upperod}{\longleftarrow}   \hskip 2mm  \uparrow^{2m-1} \oH \, ,
\]
where $\od (h) =(-1)^{\deg (h)} (h\ot 1)\cdot \oD$. We infer that,
additively,
\begin{equation}
\label{eq:circ2}
H^*(F(\oX ,2);\bbC ) =\oH^{\ot 2}/(\oD) \hskip 2mm \oplus \uparrow^{2m-1}
\Ann(\oH^+)\, ,
\end{equation}
where $\Ann(\oH^+):= \{ h\in \oH \, ;\, hh'=0, \, \forall h'\in \oH ,\
\deg (h')>0\}$. The full multiplicative structure of
$H^*F(\oX ,2)$ can in fact also be described, using
decomposition~(\ref{eq:circ2}).
\mbox{  }\hfill \endex
}
\end{example}

\begin{example}
\label{ex:circ3}
{\rm
The underlying chain complex of our punctured model $E_3(\oH)$ for $3$
points takes the  shape which we already know from Example~\ref{full},
namely

\begin{center}
{
\unitlength=1pt
\begin{picture}(180,130)
\put(140.00,3.00){\makebox(0.00,0.00){\scriptsize $\ssoD_{21}$}}
\put(140.00,55.00){\makebox(0.00,0.00){\scriptsize $-\ssoD_{32}$}}
\put(140.00,85.00){\makebox(0.00,0.00){\scriptsize $-\ssoD_{31}$}}
\put(140.00,118.00){\makebox(0.00,0.00){\scriptsize $\ssoD_{21}$}}
\put(30.00,20.00){\makebox(0.00,0.00){\scriptsize $\ssoD_{32}$}}
\put(40.00,70.00){\makebox(0.00,0.00){\scriptsize $\ssoD_{21}$}}
\put(30.00,105.00){\makebox(0.00,0.00){\scriptsize $\ssoD_{31}$}}
\put(200.00,30.00){\makebox(0.00,0.00){$\oH \cdot G_{21} G_{32}$}}
\put(200.00,90.00){\makebox(0.00,0.00){$\oH \cdot G_{21} G_{31}$}}
\put(160.00,20.00){\vector(-3,-1){52}}
\put(160.00,40.00){\vector(-3,1){52}}
\put(160.00,80.00){\vector(-3,-1){52}}
\put(160.00,100.00){\vector(-3,1){52}}
\put(80.00,0.00){\makebox(0.00,0.00){$\otexp \oH2 \cdot G_{32}$}}
\put(80.00,60.00){\makebox(0.00,0.00){$\otexp \oH2 \cdot G_{21}$}}
\put(80.00,120.00){\makebox(0.00,0.00){$\otexp \oH2 \cdot G_{31}$}}
\put(52,8){\vector(-1,1){42}}
\put(52,60.00){\vector(-1,0){35}}
\put(52,112){\vector(-1,-1){42}}
\put(0.00,60.00){\makebox(0.00,0.00){$\otexp \oH3$}}
\end{picture}}
\end{center}
or, in a more condensed form,
\[
\oH^{\ot 3} \stackrel{\od}{\longleftarrow} \oH^{\ot 2}\cdot
G_{21} \oplus \oH^{\ot 2}\cdot G_{31} \oplus \oH^{\ot 2}\cdot
G_{32} \stackrel{\od}{\longleftarrow} \oH \cdot G_{21}G_{31}
\oplus \oH \cdot G_{21}G_{32}\, ,
\]
where
\begin{eqnarray}
\label{eli5}
\lefteqn{\od \left(
(h^3_1 \ot 1\ot h^3_3) G_{21} +
(h^2_1 \ot h^2_2 \ot 1) G_{31} +
(h^1_1 \ot h^1_2 \ot 1) G_{32}
\right)=}
\\
\nonumber
&&=
(-1)^{\deg (h^3_1)+\deg (h^3_3)} (h^3_1 \ot 1\ot h^3_3) \oD_{21}+
(-1)^{\deg (h^2_1)+\deg (h^2_2)} (h^2_1 \ot h^2_2 \ot 1) \oD_{31}+
\\
\nonumber
&&
\hskip 1mm +(-1)^{\deg (h^1_1) +\deg (h^1_2)} (h^1_1 \ot h^1_2 \ot 1) \oD_{32}
\end{eqnarray}
and
\begin{eqnarray}
\label{eli4}
\lefteqn{
\od (h'\cdot G_{21}G_{31}+ h''\cdot G_{21}G_{32})=}
\\
\nonumber
&&=
-\left(
(-1)^{\deg (h')} (h'\ot 1 \ot 1) \oD_{31} +
(-1)^{\deg (h'')} (h''\ot 1 \ot 1) \oD_{32}
\right) G_{21} +
\\
\nonumber
&&
\hskip 5mm + (-1)^{\deg (h')} (h'\ot 1 \ot 1) \oD_{21} G_{31} +
(-1)^{\deg (h'')} (h'' \ot 1 \ot 1) \oD_{21} G_{32}.
\end{eqnarray}

Let us apply the above three-points punctured model to a $1$-connected smooth
compact complex algebraic surface $X$.
Note that $r:= b_2(X)\ge 1$, and that any $b_2 \ge 1$ may be
realized, for example by blowing up points in
$\bbC \bbP^2$; see \cite[1.1.1]{FM}. Here
computations are again easy, since $\oH = \bbC \cdot 1
\oplus H^2$, and $\oD = \sum^r_{i=1} x_i \ot x_i$,
where $\{ \Rada x1r \}$ denotes a convenient $\bbC$-basis of $H^2$.

It is straightforward
to see that the $2$-cycles of $E_3(\oH)$, with respect to exterior degree,
are given, in~(\ref{eli4}), by the conditions $h', h''\in H^2$.
Likewise, the exterior degree $1$-cycles of $E_3(\oH)$ are given,
in~(\ref{eli5}), by:
$h^3_1, h^2_1, h^1_2 \in H^2 $, if $r>1$. For $r=1$,
exceptional $1$-cycles of the form
\[
(1\ot 1 \ot h^3_3)  G_{21}+
(1\ot h^2_2 \ot 1) G_{31}+ (h^1_1 \ot 1 \ot 1)\ G_{32},
\]
where $h^1_1, h^2_2, h^3_3 \in H^2 $ and
$h^1_1+ h^2_2 +h^3_3 =0$, must be added. Finally,
 the bigraded Poincar\' e polynomial is given by:
\[
P_{F(\ssoX ,3)}(s,t)\! =\!
\cases
{\!\!(1+rt^2)[1+2rt^2+(r^2-3)t^4]+ st^5[(3r+ (3r^2-2)t^2]+ 2rs^2t^8}
{\hskip -2mm for $r>1$,}{\!\!(1+3t^2)+ st^5(5+t^2)+ 2s^2t^8}
{\hskip -2mm for $r=1$.}
\]
\endex
}
\end{example}

\begin{remark}
\label{xox}
{\rm
Even though the natural projection, $ F(\oX ,3)\to \oX $, always has a
section (see~\cite[Lemma~II\,1.1]{FH}), the associated Serre spectral sequence
need not collapse. Indeed, for $ X=\bbC \bbP^2 $, the Poincar\' e
polynomial \hskip 1mm $P_{\lloX }(t)=1+t^2$
\hskip 1mm does not divide $P_{F(\lloX ,3)}(t) $;
see Example~\ref{ex:circ3}.
\endex
}
\end{remark}

\section{\bf Proofs related to the $J$-model}
\label{proofs}

Let $H$ be an even-dimensional Poincar\' e duality algebra over a field $\bbK$.
In this section we prove Propositions~\ref{eli3} and~\ref{oH}, and
Theorem~\ref{main}. Let us start with Proposition~\ref{eli3}, whose proof
is the easiest.

\begin{proof}[Proof of Proposition~\ref{eli3}]
By construction, the differential tensor product $ \bbK \ot_{(H,0)}(J_n(H),d)$
is obtained from the DGA $ (E_n(H),d) $, by first killing $ G_{i1} $ and
$ dG_{i1}=\Delta _{i1} $, for $ n\geq i>1 $, and then also killing
$ \iota _1(h) $, for $ h\in H^+ $. It is now easy to check that one obtains
in this way a DGA isomorphic to $ (E_{n-1}(\oH),\od ) $.
\end{proof}

\begin{proof}[Proof of Proposition~\ref{oH}]
Let us show first that the composition
\begin{equation}
  \label{notebook}
H \ot \oHsuper {\ot l-1} \stackrel{i}{\hookrightarrow} \otexp Hl
\stackrel\pi\longrightarrow \otexp Hl/(\fatDelta ^l)
\end{equation}
is an {\em epimorphism\/}.
To this end, define the filtration
$\{F_k:=\otexp Hk \ot \oHsuper{\ot l-k}\}_{1 \leq k \leq l}$
of $\otexp Hl$ and prove, by induction
on $k$, that $\Im(\pi i)$ contains $\pi(F_k)$ for each $1 \leq k \leq
l$. Because $F_l = \otexp
Hl$, this would imply the statement.

Because $F_1 = \Im(i)$, $\Im(\pi i)$ contains $\pi(F_1)$. Suppose
that we have already proved that $\Im(\pi i) \supset
\pi(F_{k-1})$, for some $2 \leq k < l$, and let $h\in F_k\setminus
F_{k-1}$. We may clearly assume, without loss of generality, that
\begin{equation}
  \label{trick1}
h = h_1 \ot \cdots \ot h_{k-1} \ot \omega \ot h_{k+1} \ot \cdots  \ot h_l.
\end{equation}
Note then that
\begin{equation}
  \label{trick2}
h' :=  h - \Delta _{k1}(h_1\ot \dots \ot 1\ot \dots \ot h_l)\in
F_{k-1}.
\end{equation}
Since by definition $\pi(h) = \pi(h')$, the induction gives $\pi(h) \in
\Im(\pi i)$, therefore  $\Im(\pi i) \supset \pi(F_k)$.

Let us prove that the composition~(\ref{notebook}) is
{\em monic\/}. Suppose therefore that $h \in \Ker(\pi i)$, that is
\[
h=  \sum h_1 \ot \dots \ot h_l \in
F_1\cap (\fatDelta ^l) = \Im(i) \cap (\fatDelta ^l)
\]
and prove that then $h$ actually equals $0$. Clearly $h$ must be of the form
\begin{equation}
\label{eq:fi}
   h =\sum_{1 < s \leq l} \Delta _{s1} \cdot
           (\sum \varphi_{1}^{(s)} \ot \dots \ot \varphi_{l}^{(s)})\, ,
\end{equation}
for some $ \Phi ^{(s)}:=\sum \varphi_{1}^{(s)} \ot \dots \ot
 \varphi_{l}^{(s)}   \in \otexp Hl $. We infer from~(\ref{trick1})
and~(\ref{trick2}) that $\Delta _{s1} \cdot F_k \subset \Delta _{s1}\cdot
F_{k-1}+\Delta _{k1}\cdot F_k $, whenever  $ l\geq k>s>1 $, thus
we can assume that $ \Phi ^{(s)}  \in F_s $ for $1<s \leq l$.

The last possibly nonzero term
of~(\ref{eq:fi}), say $ \Delta _{q1}\cdot \Phi ^{(q)} $,
equals, by property~(\ref{eli2}) of the diagonal,
\begin{equation}
  \label{trick3}
\Delta _{q1}\cdot (\sum  \epsilon \cdot \varphi _{1}^{(q)}
\varphi _{q}^{(q)}\ot \dots \ot \varphi _{q-1}^{(q)}\ot  1\ot
\varphi _{q+1}^{(q)}\dots \ot  \varphi _{l}^{(q)}),
\end{equation}
where $\epsilon$ is an appropriate sign,
and it is zero $ \mod \, F_{q-1} $, by assumption. At the same time, the sum
from~(\ref{trick3}) equals, $ \mod \, F_{q-1} $,
$\sum  \epsilon\cdot \varphi _{1}^{(q)}
 \varphi_{q}^{(q)}\ot \dots \ot \omega \ot \dots \ot  \varphi _{l}^{(q)} $,
since $ \Phi ^{(q)}\in F_q $.
We conclude that $\sum  \epsilon\cdot \varphi _{1}^{(q)}
\varphi _{q}^{(q)}\ot \dots \ot 1\ot \dots \ot  \varphi _{l}^{(q)}=0 $,
which implies  $ \Delta _{q1}\cdot \Phi ^{(q)}=0 $.
The proposition is proved.
\end{proof}

\begin{proof}[Proof of Theorem~\ref{main}]
We will work with the \kriz-model presented as the direct sum
in~(\ref{22}). The main trick is to write $E_n(H)$ as a bicomplex. Put
\[
 E^n_{pq} := \bigoplus \otexp H{n-(p+q)}\cdot
G_{i_1j_1}G_{i_2j_2}\cdots
G_{i_{p+q}j_{p+q}},
\]
with the sum running over $\rada{i_1,j_1}{i_{p+q},j_{p+q}} \in
I_{p+q}$ as in~(\ref{22}) such that $\mbox{\it card}\{s;\ j_s = 1\} =
q$.
Then
\[
E_n(H) = \bigoplus_{0 \leq p+q \leq n-1} E^n_{pq}
\]
and the differential $d$ clearly splits as $d = d_1 + d_2$, where
\[
d_1 : E^n_{pq} \to E^n_{p,q-1}
\mbox { and }
d_2 : E^n_{pq} \to E^n_{p-1,q}.
\]
Let us agree for the rest of this proof to understand by degree
the {\em bicomplex degree,\/} that is, elements of $E^n_{pq}$
have degree $p+q$.

\begin{example}
\label{presto}
{\rm
The \kriz-model $E_3(H)$ for three points, described explicitly in
Example~\ref{full}, can be organized into a bicomplex whose nontrivial
part is:
\begin{center}
\unitlength=.4cm
\begin{picture}(10,11)(0,0)
\put(0,0){\makebox(0.00,0.00){$\otexp H3$}}
\put(10,0){\makebox(0.00,0.00){$\otexp H2\cdot G_{32}$}}
\put(0,5){\makebox(0.00,0.00){$\otexp H2\cdot G_{21}
                               \oplus \otexp H2\cdot G_{31}$}}
\put(10,5){\makebox(0.00,0.00){$H\cdot  G_{21}G_{32}$}}
\put(0,10){\makebox(0.00,0.00){$H\cdot  G_{21}G_{31}$}}
\put(10,10){\makebox(0.00,0.00){$0$}}
\put(7,0){\vector(-1,0){5}}
\put(7,5){\vector(-1,0){1.8}}
\put(8,10){\vector(-1,0){5}}
\put(0,4){\vector(0,-1){3}}
\put(0,9){\vector(0,-1){3}}
\put(10,4){\vector(0,-1){3}}
\put(10,9){\vector(0,-1){3}}
\put(.4,2.5){\makebox(0.00,0.00)[l]{\mbox {\small $d_1$}}}
\put(.4,7.5){\makebox(0.00,0.00)[l]{\mbox {\small $d_1$}}}
\put(10.4,2.5){\makebox(0.00,0.00)[l]{\mbox {\small $d_1$}}}
\put(10.4,7.5){\makebox(0.00,0.00)[l]{\mbox {\small $d_1$}}}
\put(4.5,.4){\makebox(0.00,0.00)[b]{\mbox {\small $d_2$}}}
\put(6.2,5.4){\makebox(0.00,0.00)[b]{\mbox {\small $d_2$}}}
\put(5.5,10.4){\makebox(0.00,0.00)[b]{\mbox {\small $d_2$}}}
\put(25.1,0){\makebox(0.00,0.00){\rule{10pt}{3pt}}}
\end{picture}
\end{center}
}
\end{example}

Theorem~\ref{main} will follow from the following statement.

\begin{lemma}
\label{13}
The columns of the bicomplex $(E^n_{pq},d)$ are acyclic in positive
$q$-degrees. The~only nontrivial homology of the $p$-th column is
\begin{equation}
\label{0}
H_0(E^n_{p*},d_1) \cong
\bigoplus_{J_p} \otexp H{n-p}/(\fatDelta^{n-p}) \cdot
G_{i_1j_1}G_{i_2j_2} \cdots G_{i_pj_p},
\end{equation}
where $J_p$ was introduced in Proposition~\ref{ifgod}.
\end{lemma}

Assuming Lemma~\ref{13}, Theorem~\ref{main} immediately follows from
the observation that $J_n(H)$ in presentation~(\ref{alt}) is the
second term of the obvious spectral sequence related to the bicomplex
$(E^n_{pq}, d_1 + d_2)$ and that this sequence is concentrated on the
line $ q=0 $.
It remains to prove Lemma~\ref{13}.

Since formula~(\ref{0}) for the $0$-th homology is obvious, we need
only
to prove the acyclicity in positive degrees. This will be done in two
steps.

{\bf Step~(1)} Let $T^n_*$ be the extreme left column $E^n_{0*}$ of
$E^n_{**}$. We show that all remaining columns $E^n_{p*}$, $p \geq 1$,
are combinations of complexes $T^{s}_*$ with $s < n$. So it is enough
to prove only the acyclicity of $T^n_*$, for all $n \geq 1$.

{\bf Step~(2)} To prove the acyclicity of $T^n_*$, we observe that
$T^n_*$
decomposes into the direct sum of two copies of $T^{n-1}_*$. Using
this we reduce the proof of the acyclicity of $T^n_*$ to the
verification
that a certain very explicit map is monic.

Let us start with Step~(1). For $n \geq 1$ we denote
\[
(T^n_*,d) := (E^{n}_{0*},d_1).
\]

\begin{claim}
\label{c1}
For any $p \geq 1$, the column $(E^n_{p*},d_1)$ decomposes as
\[
(E^n_{p*},d_1)
= \bigoplus_{J_p}
(E^{n-p}_{0*} \cdot G_{i_1j_1}G_{i_2j_2} \cdots G_{i_pj_p},d_1)
\cong \hskip 2mm
\uparrow^{p(2m-1)} \bigoplus_{J_p}(T^{n-p}_*,d),
\]
where $J_p$ was defined in Proposition~\ref{ifgod}.
\end{claim}

The claim is obvious, because the differential $d_1$ by definition
does not affect generators $G_{ij}$ with $j \geq 2$.

\begin{example}
{\rm The complex $(T^1_*,d)$ is just $H$ with trivial differential.
The complex $(T^2_*,d)$~is
\[
\otexp H2 \stackrel{d}{\longleftarrow} H \cdot G_{21},
\]
with $d$ given by the multiplication with $\Delta_{21}$. We see that
the right column of the bicomplex $E^3_{**}$ described in
Example~\ref{presto} is isomorphic to
\[
(T^2_*,d) \cdot G_{32} \cong \hskip 2mm\uparrow^{2m-1} (T^2_*,d),
\]
as predicted by Claim~\ref{c1}.
}
\end{example}

Let us move to Step~(2), that is, prove that the complexes
$(T^n_*,d)$, $n \geq 1$, are acyclic in positive degrees.
It follows from the definition that, for $0 \leq k \leq n-1$,
\[
T^n_k = \bigoplus \otexp H{n-k}\cdot G_{i_11}G_{i_21}\cdots G_{i_k1},
\]
where the summation runs over all $\Rada i1k$ with
$2 \leq i_1 < i_2 \cdots < i_k \leq n$.

Clearly $T^n_*$ splits into a two-column bicomplex $T^n_{**} \cong
T^n_{0*} \oplus T^n_{1*}$, with $T^n_{0q}$ consisting of summands
\[
\otexp H{n-q}\cdot G_{i_11}G_{i_21}\cdots G_{i_q1}
\]
with $i_q < n$
and $T^n_{1q}$ consisting of summands
\[
\otexp H{n-(q+1)}\cdot G_{i_11}G_{i_21}\cdots G_{i_{q+1}1}
\]
of $T^n_{1+q}$ with $i_{q+1} =
n$. The differential $d$ obviously decomposes as $d = d_1 +d_2$, where
\[
d_1 : T^n_{pq} \to T^n_{p,q-1}
\mbox { and }
d_2 : T^n_{pq} \to T^n_{p-1,q}.
\]
The following claim is evident, see also Example~\ref{jsem}.

\begin{claim}
\label{tak}
One has the following isomorphisms of complexes:
\[
(T^n_{0*},d_1) \cong (T_*^{n-1},d)\ot (H,d=0)
\mbox { and }
(T^n_{1*},d_1) \cong \hskip 2mm \uparrow^{2m-1}(T^{n-1}_*,d).
\]
\end{claim}

\begin{example}
{\rm
\label{jsem}
The complex $(T^4_*,d)$ splits as:
\begin{center}
\unitlength=.4cm
\begin{picture}(14,11)(0,-.3)
\put(0,0){\makebox(0.00,0.00){$\otexp H4$}}
\put(14,0){\makebox(0.00,0.00){$\otexp H3\cdot G_{41}$.}}
\put(0,5){\makebox(0.00,0.00){$\otexp H3\cdot G_{21}
                               \oplus \otexp H3\cdot G_{31}$}}
\put(14,5){\makebox(0.00,0.00){$\otexp H2\cdot G_{21}G_{41}
                               \oplus \otexp H2\cdot G_{31}G_{41}$}}
\put(0,10){\makebox(0.00,0.00){$\otexp H2\cdot  G_{21}G_{31}$}}
\put(14,10){\makebox(0.00,0.00){$H \cdot G_{21}G_{31}G_{41}$}}
\put(11,0){\vector(-1,0){9.5}}
\put(7,5){\vector(-1,0){1.8}}
\put(10,10){\vector(-1,0){6.5}}
\put(0,4){\vector(0,-1){3}}
\put(0,9){\vector(0,-1){3}}
\put(14,4){\vector(0,-1){3}}
\put(14,9){\vector(0,-1){3}}
\put(.4,2.5){\makebox(0.00,0.00)[l]{\mbox {\small $d_1$}}}
\put(.4,7.5){\makebox(0.00,0.00)[l]{\mbox {\small $d_1$}}}
\put(14.4,2.5){\makebox(0.00,0.00)[l]{\mbox {\small $d_1$}}}
\put(14.4,7.5){\makebox(0.00,0.00)[l]{\mbox {\small $d_1$}}}
\put(7,.4){\makebox(0.00,0.00)[b]{\mbox {\small $d_2$}}}
\put(6.4,5.4){\makebox(0.00,0.00)[b]{\mbox {\small $d_2$}}}
\put(7,10.4){\makebox(0.00,0.00)[b]{\mbox {\small $d_2$}}}
\end{picture}
\end{center}

We recognize the left column as the left column of the bicomplex in
Example~\ref{presto} tensored with $H$, that is, $T^3_* \ot
H$. The right column is $T^3_* \cdot G_{41}$, that is,
$\uparrow^{2m-1} T^3_*$. These observations are in a perfect harmony
with
Claim%
~\ref{tak}.
\endex
}
\end{example}

Assuming inductively that $(T^{n-1}_*,d)$ is acyclic in positive
dimensions, the cohomology of $(T^{n}_*,d)$ reduces to the cohomology
of a two-term complex
\[
H_0(T^{n-1}_*,d)\ot H
\longleftarrow \hskip 2mm \uparrow^{2m-1} H_0(T^{n-1}_*,d)
\]
with the differential induced by $d_2$. The above map can be easily
identified with
\[
\otexp H{n-1} /(\fatDelta^{n-1})\ot H  \longleftarrow \hskip 2mm
\uparrow^{2m-1}
\otexp H{n-1}/(\fatDelta^{n-1})
\]
where the differential is given by the multiplication with
$\Delta_{n1}$.
So it remains to prove only that the above two-term complex is acyclic
at the right term, that is

\begin{claim}
The map
\[
\otexp H{n-1}/(\fatDelta^{n-1}) \longrightarrow \otexp
H{n-1}/(\fatDelta^{n-1})\ot H
\]
sending $[\alpha] \in \otexp H{n-1}/(\fatDelta^{n-1})$ to
$(\pi \ot id)((\alpha \ot 1 )\Delta_{n1})\in \otexp H{n-1}/(\fatDelta^{n-1})
\ot H $ is a monomorphism.
\end{claim}

The proof of the claim is simple.
Suppose that $\alpha$ is homogeneous of degree~$s$.
The only term in $(\alpha \ot 1)\Delta_{n1}$ of bidegree $(s,2m)$ is
$ \alpha \ot \omega $, where $\omega \in H^{2m}$ is the fundamental
class. Therefore $(\pi \ot id)((\alpha \ot 1 )\Delta_{n1})=0 $
if and only if $[\alpha] = 0$.
This finishes the proof of Lemma~\ref{13} and thus also the proof of
Theorem~\ref{main}.
\end{proof}

\section{\bf Configuration spaces of punctured manifolds}
\label{sec:punct}

\noindent This section is devoted to the proof of Theorem~\ref{thm:punct}.
The bridge between topology and algebra is provided by the natural fibration
\begin{equation}
\label{eq:fibr}
F(\oX ,n-1)\hookrightarrow  F(X,n)\stackrel p\longrightarrow  X \ni \pt,
\end{equation}
where $ p $ is the projection onto the first coordinate. When $ X $
is as in our theorem,~(\ref{eq:fibr}) above is a fibration of connected
spaces, having a 1-connected base with finite Betti numbers. Therefore,
the approach initiated by J.C.~Moore, to get cochain algebra
models for the fiber
via differential homological algebra, may be used, in the particularly
convenient form of the so-called {\em semifree resolutions} introduced by
F\' elix-Halperin-Thomas in~\cite{RHT}.

Parts of the theory used
below work for arbitrary field coefficients (e.g., the connection
between semifree resolutions and fibrations), some other parts
(e.g., the theory of Sullivan relative models for fibrations)
require characteristic zero field coefficients, while Morgan's theory
of bigraded models for smooth complex varieties from~\cite[Section 9]{M}
needs $ \bbC $-coefficients. To simplify matters, we will thus work over
$ \bbC $, throughout this section.

\subsection{Semifree resolutions, relative models, and fibrations}
\label{subsec:tools}
We start by recalling some relevant facts from~\cite{RHT}.
{}From now on, all DGA's $ (A,d_A) $ will be tacitly assumed to be
homologically connected, that is, $ H^0(A,d_A)=\bbC. $ Given a DGA map,
$ (B,d_B)\stackrel f \longrightarrow (E,d_E) $, we will assume
in \S~\ref{subsec:tools} that
$ H^1f $ is monic. The examples we have in mind are
\begin{equation}
\label{eq:derham}
\dR (X)\stackrel {\dR (p)}\longrightarrow \dR (F(X,n)) ,
\end{equation}
coming from the basic fibration (\ref{eq:fibr}), via the
$ C^{\infty}$ de Rham DGA functor with complex coefficients, $ \dR
(\cdot )$. The following definition is taken from~\cite[\S~6]{RHT}.

\begin{definition}
\label{def:sfree}
A DGA $ (\BB ,d_{\BB}) $ is called {\em semifree\/} over $(B,d_B)$ if
there exists a graded vector space $Z^*$ with the following
properties:
\begin{itemize}
\item[(i)]
$\BB ^*=B^*\otimes Z^* $, as graded vector spaces,
\item[(ii)]
the graded vector space $Z^*$ admits a second grading,
$Z^*=\bigoplus _{k\geq 0}Z^*_k $,
such that $ Z^*_0=\bbC $ and $ (B\otimes Z_0,d_{\BB})=(B,d_B) $ as
DGA's,
\item[(iii)]
the product of $\BB$ satisfies $ b\ot z=(b\ot 1)\cdot
(1\ot z) $, for $ b\in B $ and $ z\in Z $,
\item[(iv)]
all subspaces $ B\otimes Z_{\leq k} $ are $ d_{\BB}$-invariant, and, finally
\item[(v)]
for $k>0$, the  quotient cochain complexes
are of the form
\[
(B\otimes Z_{\leq k}/B\otimes Z_{<k}, d_{\BB})=(B,d_B)\otimes (Z_k,0)\, .
\]
\end{itemize}
A DGA map
$
(\BB,d_{\BB})\stackrel {\rho }\longrightarrow (E,d_E),
$
is said to be a {\em semifree resolution\/} of $ f: (B,d_B)\to (E,d_E)
$ if $ (\BB ,d_{\BB}) $ is semifree over $ (B,d_B) $, the restriction
of $ \rho $ to $ (B\otimes Z_0, d_{\BB})\equiv (B,d_B) $ equals $f$,
and $ H^*\rho $ is an isomorphism.
\end{definition}

The main source of semifree resolutions is provided by {\em Sullivan
relative models\/} of $ f $. That is, by DGA's of the form
$ (\BB,d_{\BB})=(B\otimes \ext W^*,d_{\BB}) $, where
$ W^*=\bigoplus _{k>0}W^k $,
and the differential $ d_{\BB} $ satisfies
a certain nilpotence condition, together with a DGA map,
$ (\BB,d_{\BB})\stackrel {\rho }\longrightarrow (E,d_E) $, which
restricts to $ f $ on $ B $ and induces a homology isomorphism.
If $ H^1f $ is monic, then $ f $ has such a relative model, which
is a semifree resolution, in the sense of Definition~\ref{def:sfree}.
See~\cite[Proposition 14.3 and Lemma 14.1]{RHT}.

 Assume now that $ B $ is
{\em augmented}, by a DGA map $ \varepsilon :(B,d_B)\rightarrow \bbC $.
In the geometric case from (\ref{eq:derham}),
$B^* = \dR (X)$ and the augmentation comes
from the inclusion, $ \pt \hookrightarrow X $, via $ \dR (\cdot ) $.
Moore's basic result which relates topology to algebra says that the
DGA
\begin{equation}
\label{eq:moore}
(\ext W,\overline d):=\bbC\otimes _{(B,d_B)}(\BB,d_{\BB})
\end{equation}
is a Sullivan model of the fiber, $ F(\oX,n-1) $, for any relative
model, $ (\BB,d_{\BB}) $, of the DGA map $ \dR (p) $ from
(\ref{eq:derham}). Here, $ \bbC $ is to be considered as a right
$ (B,d_B) $-module, via $ \varepsilon $, and $ (\BB,d_{\BB}) $ as a
left $ (B,d_B) $-module, via the inclusion $ (B,d_B)\hookrightarrow
(\BB,d_{\BB}) $. See~\cite[Theorem 15.3 and Theorem 7.10]{RHT}.
Note also that
this is the only place where we need the 1-connectivity assumption on
$ X $.

\subsection{\bf Algebraic weak homotopy type and homotopy fiber}
\label{subsec:weak}
Our strategy for proving Theorem~\ref{thm:punct} is to use the
differential tensor product construction from (\ref{eq:moore}) to
arrive at the DGA model of $ F(\oX,n-1) $ from our theorem. We are
thus led to study how this construction depends on its input,
that is, on the DGA map $ \dR (p) $.
The next definition comes in  naturally.

\begin{definition}
\label{def:htype}
Let $ f:B \to E $ and $ f':B'\to E' $ be DGA maps,
with $ B $ and $ B' $ augmented. An {\em elementary
weak equivalence\/}, $ f\to f' $, consists of a DGA map,
$ \Phi :E\to E' $, and an augmented DGA map, $ \varphi :B\to B'$,
such that both $ H^*\Phi $ and $ H^* \varphi $ are isomorphisms,
and such that $ \Phi f $ and $ f'\varphi $ are Sullivan homotopic
DGA maps (notation: $ \Phi f \simeq f'\varphi $):
\begin{center}
\unitlength=.4cm
\begin{picture}(10,10.3)(-1,-4)
\put(0,0){\makebox(0.00,0.00){$ B $}}
\put(10,0){\makebox(0.00,0.00){$ B' $}}
\put(0,5){\makebox(0.00,0.00){$ E $}}
\put(10,5){\makebox(0.00,0.00){$ E' $}}
\put(5,-3.5){\makebox(0.00,0.00){$ \bbC $}}
\put(5,2.5){\makebox(0.00,0.00){$ \HH $}}
\put(4.5,6){\makebox(0.00,0.00)[l]{$ \Phi $}}
\put(4.5,-1){\makebox(0.00,0.00)[l]{$ \varphi $}}
\put(2,-2.5){\makebox(0.00,0.00)[l]{$ \varepsilon $}}
\put(7.5,-2.5){\makebox(0.00,0.00)[l]{$ \varepsilon ' $}}
\put(-0.5,2.5){\makebox(0.00,0.00)[r]{$ f $}}
\put(10.5,2.5){\makebox(0.00,0.00)[l]{$ f' $}}
\put(0,1){\vector(0,1){3}}
\put(10,1){\vector(0,1){3}}
\put(1,0){\vector(1,0){8}}
\put(1,5){\vector(1,0){8}}
\put(1,-1){\vector(4,-3){3}}
\put(9,-1){\vector(-4,-3){3}}
\end{picture}
\end{center}
We say that
$ f $ and $ f' $ {\em have the same weak homotopy type }
(notation: $ f\cong f' $) if there is a chain of DGA maps,
$ \{f_i:B_i\to E_i \}_{0\leq i \leq l} $, as above, together
with elementary weak equivalences, $ f_i \rightarrow f_{i+1} $ or
$ f_i \leftarrow f_{i+1} $, such that $ f_0=f $ and $ f_l=f' $.
\end{definition}

\begin{lemma}
\label{lem:ww}
Assume $ f\cong f'$, in the sense of the above definition, where
$ H^1f $ (and consequently also $ H^1f' $) are monic. Let
$ \rho :\BB =(B\otimes \ext W,d_{\BB})\rightarrow E $ and
$ \rho ' :\BB' =(B'\otimes \ext W',d_{\BB'})\rightarrow E' $
be arbitrary relative models, for $ f $ and $ f' $ respectively.
Then the DGA's $\bbC \otimes_{(B,d_B)} (\BB,d_{\mathcal B})$
and $\bbC \otimes_{(B',d_{B'})} (\BB',d_{\mathcal B'})$
have the same Sullivan minimal model.
\end{lemma}

\begin{proof}
Plainly, it is enough to treat the case of an elementary
weak equivalence, say, $ f \to f' $. Suppose first that it is
{\em strict }, that is, $ \Phi f=f'\varphi $ (not only equality
up to algebraic homotopy). Consider the {\em pushout\/} construction
from ~\cite[\S~14(a)]{RHT}. It provides a DGA map,
$ \varphi \otimes id :\BB=(B\otimes \wedge W,d_{\BB})\rightarrow
(B'\otimes \wedge W,d'')=:\varphi _*\BB $, extending
$ \varphi $ and inducing a homology isomorphism;
see~\cite[Lemma 14.2]{RHT}.

Define a DGA map, $ \varphi _*\rho :\varphi _*\BB \rightarrow E' $,
which extends $ f' $, by setting $ \varphi _*\rho :=\Phi \rho $,
on $ W $. By construction, $ \varphi _*\rho $ is a relative model
of $ f' $ (see~\cite[\S~14(a)]{RHT}). By~\cite[Theorem 6.10(ii)]{RHT},
$ \bbC\otimes _B \BB $ and $ \bbC\otimes _{B'}\varphi _* \BB $ have
the same minimal model. At the same time (\cite[Proposition 14.6]{RHT}),
there is a DGA map, $ \Psi :\BB' \rightarrow \varphi _*\BB $,
extending the identity on $ B' $ and such that $ \varphi _*\rho
\circ \Psi \simeq \rho' $. It follows that $ \Psi $ induces
an isomorphism on homology; hence, by~\cite[Proposition 6.7(ii)]{RHT},
$ \bbC \otimes _{B'}\BB' $ and $ \bbC \otimes _{B'}\varphi _*\BB $
have the same minimal model. This settles the strict case.

In the general case, the property $ \Phi f \simeq f'\varphi $ means
the existence of a DGA map, $ \HH :B\rightarrow E'
\otimes \ext (t,dt) $, such that $ e_0\HH =\Phi f $
and $ e_1\HH =f'\varphi $, where $ e_i:E'\otimes \ext
(t,dt)\rightarrow E' $ are the DGA maps extending
$ id_{E'} $, which send $ t $ to $ i $ and $ dt $ to $ 0 $, for $
i=0,1 $.

Pick a relative model of $ \HH $, $ \eta :\CC=(B\otimes \ext T,
d_{\CC})\rightarrow E'\otimes \ext (t,dt) $. Since both
$ e_0 $ and $ e_1 $ induce homology isomorphisms
~\cite[Lemma 12.5]{RHT},
$ e_0 \eta $ is a relative model of
$ \Phi f $, and $ e_1 \eta $ is a relative model of
$ f'\varphi $. We may now use the strict case to infer that
$ \bbC \otimes _B \BB $, $ \bbC \otimes _{B'} \BB' $ and
$ \bbC \otimes _B \CC $ have the same minimal model.
\end{proof}

We may rephrase Lemma~\ref{lem:ww} in the following way. Given a
DGA map, $ f:B\to E $, with $ H^1f $ monic, we may define its
{\em homotopy fiber}, $ \FF(f) $, to be the minimal model of
$ \bbC \otimes _B \BB $, where $ \BB $ is any relative model of
$ f $. If $ f \cong f' $, then $ \FF (f)=\FF (f') $. In
particular, if $ \dR (p) $ from (\ref{eq:derham}) has the
same weak homotopy type as $ f $, then $ \bbC \otimes _B \BB $ will
be a (not necessarily minimal) model of $ F(\oX ,n-1) $, according
to our discussion from the end of \S~\ref{subsec:tools}. Recall that $H$
denotes the cohomology algebra $ H^*(X;\bbC ) $.

\begin{lemma}
\label{lem:red}
Let $ \iota _1:(E_1(H),d)=(H,0)\to (E_n(H),d) $
be the DGA map induced by the map
$\iota_1 : H\to \otexp Hn$ described in Section 1,
canonically augmented by $ \varepsilon (h)=0 $,
for $ h\in H^+ $, and $ \varepsilon (1)=1 $.
Then $ \FF(\iota _1) $ is the minimal model of $ (E_{n-1}(\oH ),\od ) $.
\end{lemma}

\begin{proof}
According to Theorem \ref{main} and Lemma \ref{lem:ww},
$ \FF (\iota _1)=\FF (\Psi \circ \iota _1) $. We claim now that
$ (J_n(H),d) $ is $ (H,0) $-semifree. Indeed, Propositions \ref{ifgod} and
\ref{oH} together imply that we may take
\[
Z_{k+1} :=\bigoplus _{J_k}
\oH ^{\otimes n-k-1}\cdot G_{i_1j_1}\cdots G_{i_kj_k},
\]
for $  1\leq k\leq n-2 $,
$Z_1 :=(\oH ^{\otimes n-1})^+$   and $Z_0 :=\bbC \cdot 1$.
To check the last semifreeness condition
from Definition ~\ref{def:sfree}, it suffices to recall from
Remark~\ref{bid} that $ d $ is homogeneous of degree
$ -1 $, with respect to exterior degree.

Let $\rho: \BB \to J_n(H)$
be a relative model of $ \Psi \circ \iota _1 $. It follows
that both $ \rho :\BB \to J_n(H) $ and $ id_{J_n(H)} $ are semifree
resolutions of $ \Psi \circ \iota _1 $. Therefore, $ \bbC \otimes _H\BB $
and $ \bbC \otimes _H J_n(H) $ have the same minimal model, by
\cite[Proposition~6.7(ii)]{RHT}. To finish the proof of our lemma, it is enough
to recall from Proposition~\ref{eli3} that the DGA  $ \bbC \otimes _H J_n(H) $
is isomorphic to $ E_{n-1}(\oH )$.
\end{proof}

\subsection{Replacing forms by $ E $-models}
\label{subsec:hodge}

Let $X$ be a smooth compact, complex algebraic variety. Set $ H=H^*(X;\bbC ) $.
We know, from~\cite{M} and~\cite{FMP, kriz} respectively, that
$ (H,0) $ models  $ \dR (X)$ (the formality property),
and $ (E_n(H),d) $ models the de Rham DG-algebra
$ \dR (F(X,n)) $. The next proposition links these two things,
into the form of a statement about DGA maps. By Lemma~\ref{lem:red}
and the discussion preceding it, this statement implies
Theorem~\ref{thm:punct}.

\begin{proposition}
\label{lem:rel}
The DGA map $ \dR (p) :\dR (X)\to \dR (F(X,n))  $
has the same weak homotopy type as $ \iota _1:E_1(H)\to E_n(H) $.
\end{proposition}

The main step in the proof uses  Morgan's results from~\cite{M},
so we recall his basic constructions. For $ Y=\wtY \setminus D $, a complement
of a divisor with normal crossings in a smooth compact complex algebraic
variety, Morgan introduced a {\em filtered model}, $ \MM _Y $ (denoted by
$E_{\cinf}(Y)_{\bbC}$
in~\cite[\S 2]{M}), a DG-algebra of global sections of a sheaf related
to $D$. If $ Y=\wtY $,  $ \MM _Y=\dR (Y) $. In general, there is a
DGA map,
$$ \dR (Y)\stackrel {\Psi _Y}\longleftarrow \MM _Y \, , $$
inducing an isomorphism in homology; see~\cite[\S~\S~2-3]{M}.

Moreover, $ \MM _Y $ is provided with an increasing filtration $W_*$ induced
by the stratification of $D$, and also with the associated increasing
filtration, $ Dec\, W_* $, defined in~\cite[\S 1]{M}.
Morgan's second construction is the
{\em bigraded minimal model}, $ \NN _Y $, of the filtered model $ \MM _Y $.
There is a DGA map,
$$ \NN _Y=\bigoplus_{n\geq 0}\NN ^n\stackrel {\Phi _Y }\longrightarrow
              \MM _Y \, ,$$
inducing an isomorphism in homology. Moreover,
each homogeneous component is bigraded: $ \NN ^n=\bigoplus _{r,s\in \bbN }
\NN ^{n;r,s} $, $ \NN ^0=\NN ^{0;0,0}=\bbC \cdot 1 $; the differential and the
multiplication are homogeneous, of bidegree $ (0,0) $, with respect to the
above extra bigrading of $ \NN _Y $. Finally,
 the bigrading of $ \NN _Y $  and the filtration of $ \MM _Y $ are
related by:
\begin{equation}
\label{dec}
\Phi (\NN ^{n;r,s})\subset Dec\, W_{r+s}(\MM ^n)
\end{equation}
(see~\cite[\S~6]{M}, in particular $(6.0)$ and Theorem~6.6).

Using {\em mixed Hodge diagrams}~\cite[Definition~3.5]{M} and
{\em mixed Hodge homotopies}~\cite[\S~6.1]{M}, Morgan
proved the naturality of these two
constructions. Given an algebraic map, $ f:(\wtY ,Y)\to (\wtZ ,Z) $,
one can construct a filtered DGA map $ \MM _f $, a DGA map
$ \NN _f $, which is homogeneous of  bidegree $ (0,0) $, with respect to the
extra bigradings on minimal models, and also homotopies $ \HH $ and $ \KK $
which fit into the following diagram:
\begin{center}
\unitlength=.4cm
\begin{picture}(30,3.5)(-5,4)
      \put(0,5){\makebox(0.00,0.00){$ \dR (Y) $}}
      \put(20,0){\makebox(0.00,0.00){$ \NN _Z $}}
      \put(10,0){\makebox(0.00,0.00){$ \MM _Z  $}}
      \put(0,0){\makebox(0.00,0.00){$ \dR (Z) $}}
      \put(20,5){\makebox(0.00,0.00){$ \NN _Y $}}
      \put(10,5){\makebox(0.00,0.00){$ \MM _Y $}}
      \put(5,-3.5){\makebox(0.00,0.00){$ \bbC $}}
      \put(15,-3.5){\makebox(0.00,0.00){$ \bbC $}}
      \put(5.5,2.5){\makebox(0.00,0.00){$ \HH $}}
      \put(15.5,2.5){\makebox(0.00,0.00){$ \KK $}}
      \put(5,6){\makebox(0.00,0.00)[l]{$ \Psi _Y $}}
      \put(5,-1){\makebox(0.00,0.00)[l]{$ \Psi _Z $}}
      \put(15,6){\makebox(0.00,0.00){$ \Phi _Y $}}
      \put(15,-1){\makebox(0.00,0.00){$ \Phi _Z $}}
      \put(20.5,2.5){\makebox(0.00,0.00)[l]{$ \NN _f $}}
      \put(10.5,2.5){\makebox(0.00,0.00)[l]{$ \MM _f $}}
      \put(-0.5,2.5){\makebox(0.00,0.00)[r]{$ \dR (f) $}}
      \put(0,1){\vector(0,1){3}}
      \put(10,1){\vector(0,1){3}}
      \put(20,1){\vector(0,1){3}}
      \put(9,0){\vector(-1,0){7}}
      \put(9,5){\vector(-1,0){7}}
      \put(19,0){\vector(-1,0){7}}
      \put(19,5){\vector(-1,0){7}}
      \put(1,-1){\vector(4,-3){3}}
      \put(9,-1){\vector(-4,-3){3}}
      \put(11,-1){\vector(4,-3){3}}
      \put(19,-1){\vector(-4,-3){3}}
\end{picture}
\end{center}
\begin{equation}
\label{htp}
\end{equation}
   \begin{center} \unitlength=.4cm  \begin{picture}(30,5)(-5,0)
   \end{picture}   \end{center}
Moreover, $ \KK :\NN _Z \to \MM _Y\otimes \ext (t,dt) $ has the
property that:

\begin{equation}
\label{dech}
 \KK (\NN _Z^{n;r,s})\subset Dec\, W_{r+s}((\MM _Y\otimes \ext
(t,dt))^n)\, ,
\end{equation}
where $ W_*(\MM _Y\otimes \ext (t,dt)):=W_*(\MM _Y )\otimes
\ext (t,dt) $. See Definitions~3.5 and~3.7, and Propositions~3.6 and~3.9
of~\cite{M}, for the left square from~(\ref{htp}); to construct the
right square from~(\ref{htp}), use Theorem~6.7 and Corollary~6.8 of~\cite{M}.

By construction, $ \MM _Y $ is a Sullivan model of $ Y $. Let us now
consider the decreasing DGA filtration of $ \MM _Y $, $ W^*:=W_{-*} $,
and the associated spectral sequence of DGA's,
$$ (\sir ^*(\MM _Y),d_r)_{r\geq 0}\Longrightarrow H^*(\MM _Y) \, . $$
A basic result in Morgan's theory~(\cite[Theorem~9.6]{M}) says that the
DGA $ (\ss (\MM _Y),d_1) $ is also a model of $Y$.
We need the following relative version of this result, where
$\ss (\MM_f)$ denotes the DGA map between $\ss$--terms induced by the
filtered DGA map $\MM_f$.

\begin{proposition}
\label{morg}
If $ f:(\wtY,Y)\to (\wtZ,Z) $ is an algebraic map, then the DGA map
$$ \MM _f:(\MM _Z,d)\to (\MM _Y,d) $$
has the same weak homotopy type as the DGA map
$$ \ss (\MM _f):(\ss (\MM _Z),d_1)\to (\ss (\MM _Y),d_1)\, . $$
\end{proposition}

\begin{proof}
Define decreasing DGA filtrations, on both $(\NN_Y, d)$ and
$(\NN_Z, d)$, by : $ W^k(\NN ^n):=\bigoplus _{r+s\leq n-k}\NN ^{n;r,s} $.
Relation~(\ref{dec}) implies that $\Phi$ is a filtered DGA map:
$\Phi (W^k(\NN ^n))\subset W^k(\MM^n)$. Obviously, the bihomogeneous DGA map
$\NN_f$ respects filtrations as well. It is equally easy to check that one
has DGA identifications,
$$ (\NN ,d) \equiv (\ss (\NN) , d_1) \, ,$$
for both $Y$ and $Z$, which give an identification
$$ \NN_f \equiv \ss (\NN_f) \, .$$

We claim that, when applying the $\ss (\cdot)$-functor to the right square from
(\ref{htp}), one gets  the commutative right square below, in the
DGA category.
\begin{center}
\unitlength=.4cm
\begin{picture}(30,11)(-1,-4)
    \put(16.5,5){\makebox(0.00,0.00){$ (\ss (\NN _Y),d_1) $}}
    \put(10,0){\makebox(0.00,0.00){$ (\NN _Z,d) $}}
    \put(0,0){\makebox(0.00,0.00){$ (\MM _Z,d) $}}
    \put(16.5,0){\makebox(0.00,0.00){$ (\ss (\NN _Z),d_1) $}}
    \put(10,5){\makebox(0.00,0.00){$ (\NN _Y,d) $}}
    \put(0,5){\makebox(0.00,0.00){$ (\MM _Y,d) $}}
    \put(5,-3.5){\makebox(0.00,0.00){$ \bbC $}}
    \put(21,-3.5){\makebox(0.00,0.00){$ \bbC $}}
    \put(5,2.5){\makebox(0.00,0.00){$ \KK $}}
    \put(4.5,6){\makebox(0.00,0.00)[l]{$ \Phi _Y $}}
    \put(4.5,-1){\makebox(0.00,0.00)[l]{$ \Phi _Z $}}
    \put(21.5,-1.3){\makebox(0.00,0.00){$ \ss (\Phi _Z ) $}}
    \put(21.5,6.3){\makebox(0.00,0.00){$ \ss (\Phi _Y ) $}}
    \put(26.2,0){\makebox(0.00,0.00){$ (\ss ( \MM _Z),d_1) $}}
    \put(26.2,5){\makebox(0.00,0.00){$ (\ss ( \MM _Y),d_1) $}}
    \put(26.5,2.5){\makebox(0.00,0.00)[l]{$ \ss (\MM _f) $}}
    \put(9.5,2.5){\makebox(0.00,0.00)[r]{$ \NN _f $}}
    \put(15.5,2.5){\makebox(0.00,0.00)[r]{$ \ss (\NN _f) $}}
    \put(-0.5,2.5){\makebox(0.00,0.00)[r]{$ \MM _f $}}
    \put(0,1){\vector(0,1){3}}
    \put(10,1){\vector(0,1){3}}
    \put(16,1){\vector(0,1){3}}
    \put(26,1){\vector(0,1){3}}
    \put(8,0){\vector(-1,0){6}}
    \put(8,5){\vector(-1,0){6}}
    \put(12,-0.5){$ \equiv $}
    \put(12,4.5){$ \equiv $}
    \put(20,0){\vector(1,0){2.5}}
    \put(20,5){\vector(1,0){2.5}}
    \put(1,-1){\vector(4,-3){3}}
    \put(9,-1){\vector(-4,-3){3}}
    \put(17,-1){\vector(4,-3){3}}
    \put(25,-1){\vector(-4,-3){3}}
 \end{picture}
\end{center}
Note also that the DGA maps $\ss (\Phi_Y)$ and  $\ss (\Phi_Z)$
induce homology isomorphisms, as follows from Morgan's proof
of~\cite[Theorem 9.6]{M}.

Granting the claim, we may quickly finish the proof of our proposition, as
follows. Directly from Definition~\ref{def:htype},
we deduce from the above diagram
elementary weak equivalences,
$$ \MM _f\leftarrow \NN _f\equiv \ss (\NN _f)\rightarrow \ss (\MM _f)\, , $$
and we are done.

Going back to our commutativity claim, let us start by noting that $ \KK $
is a filtered DGA map; see~(\ref{dech}). We thus get an induced DGA map,
\[
\spos (\KK ):(\spos (\NN _Z),d_0) \to (\spos (\MM _Y\otimes
\ext(t,dt)),d_0)\equiv  (\spos (\MM _Y),d_0)\otimes
\ext(t,dt).
\]
We thus obtain from the second square of~(\ref{htp})
a homotopy commutative DGA square (the left
diagram of the following display), and, for the next level of spectral
sequences, a strictly commutative one (the right diagram below):

\begin{center}
\unitlength=.4cm
\begin{picture}(33,8)(-1,-1)
   \put(21,5){\makebox(0.00,0.00){$ (\ss (\NN _Y),d_1) $}}
\put(10.2,0){\makebox(0.00,0.00){$ (\spos (\MM _Z),d_0) $}}
\put(0,0){\makebox(0.00,0.00){$ (\spos (\NN _Z),d_0) $}}
   \put(21,0){\makebox(0.00,0.00){$ (\ss (\NN _Z),d_1) $}}
\put(10.2,5){\makebox(0.00,0.00){$ (\spos (\MM _Y),d_0) $}}
\put(0,5){\makebox(0.00,0.00){$ (\spos (\NN _Y),d_0) $}}
\put(5,2.5){\makebox(0.00,0.00){$ \spos (\KK ) $}}
\put(3.3,6.2){\makebox(0.00,0.00)[l]{$ \spos (\Phi _Y) $}}
\put(3.3,-1.2){\makebox(0.00,0.00)[l]{$ \spos (\Phi _Z) $}}
\put(24.3,-1.2){\makebox(0.00,0.00)[l]{$ \ss (\Phi _Z ) $}}
\put(24.3,6.2){\makebox(0.00,0.00)[l]{$ \ss (\Phi _Y ) $}}
\put(31.2,0){\makebox(0.00,0.00){$ (\ss ( \MM _Z),d_1) $}}
\put(31.2,5){\makebox(0.00,0.00){$ (\ss ( \MM _Y),d_1) $}}
    \put(31.5,2.5){\makebox(0.00,0.00)[l]{$ \ss (\MM _f) $}}
    \put(10.5,2.5){\makebox(0.00,0.00)[l]{$ \spos (\MM _f) $}}
    \put(20.5,2.5){\makebox(0.00,0.00)[r]{$ \ss (\NN _f) $}}
    \put(-0.5,2.5){\makebox(0.00,0.00)[r]{$ \spos (\NN _f) $}}
     \put(0,1){\vector(0,1){3}}
     \put(10,1){\vector(0,1){3}}
     \put(21,1){\vector(0,1){3}}
     \put(31,1){\vector(0,1){3}}
  \put(3.5,0){\vector(1,0){3}}
  \put(3.5,5){\vector(1,0){3}}
  \put(24.5,0){\vector(1,0){3}}
  \put(24.5,5){\vector(1,0){3}}
\end{picture}
\end{center}
This verifies our claim and thus ends the proof of Proposition~\ref{morg}.
\end{proof}

\begin{proof}[Proof of Proposition~\ref{lem:rel}]
We put together a sequence of elementary weak equivalences coming
from~\cite{M},~\cite{FMP}, and~\cite{kriz}.

Applying~(\ref{htp}) (the left square) and Proposition~\ref{morg} to the
projection onto the first coordinate $ p:(X[n],F(X,n))\to (X,X) $
(where $ X[n] $ is the compactification in~\cite{FMP}), we obtain
a weak homotopy equivalence $ \dR (p)\cong \ss (\MM _p) $.

The last two elementary equivalences are described in the next diagram:

\begin{center}
\unitlength=.4cm
\begin{picture}(20,11)(-1,-4)
\put(20,5){\makebox(0.00,0.00){$ E_n(H) $}}
\put(10,0){\makebox(0.00,0.00){$ \AA ^*(X,1)_{\bbC } $}}
\put(20,0){\makebox(0.00,0.00){$ E_1(H) $ }}
\put(22,0){\makebox(0.00,0.00){$ . $ }}
\put(0,0){\makebox(0.00,0.00){$ \ss (\MM _X) $}}
\put(10,5){\makebox(0.00,0.00){$ \AA ^*(X,n)_{\bbC } $}}
\put(0,5){\makebox(0.00,0.00){$ \ss (\MM _{F(X,n)}) $}}
\put(5,-3.5){\makebox(0.00,0.00){$ \bbC $}}
\put(15,-3.5){\makebox(0.00,0.00){$ \bbC $}}
    \put(15,6){\makebox(0.00,0.00)[l]{$ \kk _n $}}
    \put(15,-1){\makebox(0.00,0.00)[l]{$ \kk _1 $}}
    \put(20.5,2.5){\makebox(0.00,0.00)[l]{$ \iota _1 $}}
\put(9.5,2.5){\makebox(0.00,0.00)[r]{$ \iota ^{\AA }_1  $}}
\put(-0.5,2.5){\makebox(0.00,0.00)[r]{$ \ss (\MM _p) $}}
\put(0,1){\vector(0,1){3}}
\put(10,1){\vector(0,1){3}}
     \put(20,1){\vector(0,1){3}}
\put(3,0){\line(1,0){4}}
\put(3,0.3){\line(1,0){4}}
\put(3.5,5){\vector(1,0){3.5}}
\put(5,6){\makebox(0.00,0.00)[l]{$ \mu _n $}}
\put(5,-1){\makebox(0.00,0.00)[l]{$ \mu _1 $}}
     \put(18,0){\line(-1,0){5}}
     \put(18,0.3){\line(-1,0){5}}
     \put(18,5){\vector(-1,0){5}}
     \put(1,-1){\vector(4,-3){3}}
     \put(9,-1){\vector(-4,-3){3}}
     \put(11,-1){\vector(4,-3){3}}
     \put(19,-1){\vector(-4,-3){3}}
\end{picture}
\end{center}

The DGA $ \AA ^*(X,n )_{\bbC} $ is the Fulton-MacPherson model from
~\cite[Theorem 8]{FMP} (with $\bbC$ coefficients). The DGA isomorphism
$\mu_n$ is constructed in~\cite[\S 6]{FMP}. The DGA's $ \AA ^*(X,1 )_{\bbC} $
and $\ss (\MM_X)$ are both equal to $(H, 0)$, and $\mu_1 = id$. The DGA
$ \AA ^*(X,n )_{\bbC} $ is a quotient of a certain free graded
$H^{\otimes n}$-algebra, and the DGA map $\iota_1^{\AA}$ between $\AA$-models
is induced by $\iota_1 : H \to H^{\otimes n}$, like in the case of $E$-models.
An easy analysis of the construction of $\mu_n$ shows that the above left
DGA square is commutative, and therefore provides an elementary weak
equivalence, $\ss (\MM_p)\to \iota_1^{\AA}$.

The DGA map $\kk_n$, constructed by \kriz \, in \cite[\S 3]{kriz}, induces
an isomorphism in homology; see~\cite[\S 4]{kriz}. By construction,
$\kk_n$ is $H^{\otimes n}$-linear. Consequently, the above right DGA square
is commutative, and gives the last needed weak equivalence,
$ \iota_1^{\AA} \leftarrow \iota _1 $.

Thus, Proposition~\ref{lem:rel} is proved, and the proof of
Theorem~\ref{thm:punct} is complete.
\end{proof}

\section{\bf A panorama of natural structures on punctured manifolds}
\label{sec:struct}

In this section, we will define various additional structures on
$E$-models, with emphasis on the punctured case. They seem to be
both highly natural, and potentially useful for further applications.

Our starting point is like in Section~\ref{intro}: a Poincar\' e
duality algebra, $ H $, over a field $ \bbK $,
together with an orientation class, $ \omega \in H^{2m}\setminus \{ 0\} $.
As noted in Remark~\ref{bid}, the associated models, $ E_n(H) $
and  $ E_n(\oH ) $, are bigraded differential algebras (DBGA's).

\subsection{Symmetry}
\label{subsec:sym}

The symmetric group $ \Sigma _n $ acts (on the left) on both
$ E_n(H) $ and  $ E_n(\oH ) $, by DBGA maps, and the canonical projection,
$ E_n(H)\to  E_n(\oH ) $, is a $ \Sigma _n $-equivariant DBGA map.
Geometrically, the above projection corresponds to the natural
inclusion, $ F(\oX ,n)\hookrightarrow F(X,n) $, which is equivariant
with respect to the natural $ \Sigma _n $-actions on $ n $-configurations.

Algebraically, the action of $ \sigma \in \Sigma _n $ on the algebra
generators is defined as follows. For
elements $ h_1,\dots ,h_n $ of $H$ or $ \oH $,
$ \sigma \cdot h_1\otimes \dots \otimes h_n:=\pm h_{\sigma ^{-1}1}\otimes
\dots \otimes h_{\sigma ^{-1}n} $, where the sign depends on $ \sigma $
and the degrees of $ \{ h_i \} $, according to the standard Koszul
sign convention. For both $ E_n(H) $ and  $ E_n(\oH ) $, $ \sigma \cdot
G_{ij}:=G_{\sigma i,\sigma j} $. Here, it is convenient to replace
the definitions from \S~\ref{intro} by the following equivalent ones:
take as exterior generators {\em all\/} $ \{ G_{ij}; n\geq i\neq j
\geq 1 \} $, and add the relations $ G_{ij}=G_{ji} $, for all $ i\neq j $.
As for the projection, $ E_n(H)\to  E_n(\oH ) $, it is induced by the
canonical projection, $ H\to \oH $, on $ H^{\otimes n} $,
and acts as the identity on the exterior generators, $ G_{ij} $.

\subsection{Simplicial structure}
\label{subsec:simpl}

Set $ E_0(\oH ):=\bbK $, and $ \EE (\oH ):=\{ E_n(\oH )\} _{n\geq 0} $.
We will give $ \EE (\oH ) $ the structure of a {\em simplicial DBGA }.
For $ H=H^*(S^2) $, this structure is given by Drinfel'd's \cite[p.843]{D}
natural {\em cosimplicial group } structure of classical (Artin)
pure braid groups, $\{ \pi _1F(\bbR ^2,n) \}$, via the cohomology
algebra functor; see also Example~\ref{ex:OS}.

Set $ D_0=D_1:=\varepsilon :
(\oH ,0)=E_1(\oH )\to E_0(\oH )=\bbK $, where $ \varepsilon $ denotes
the canonical augmentation of the connected graded algebra $ \oH $. Set
also $ S_0:=\eta :E_0(\oH )\to E_1(\oH ) $, where $ \eta $ is the
unit of $ \oH $.

Fix now $ n\geq 1 $. The degeneration maps, $ S_k:E_n(\oH )\to
E_{n+1}(\oH ) $, where $ 0\leq k\leq n $, correspond to the natural
projections, $ pr_{k+1}:F(\oX ,n+1)\to F(\oX ,n) $, which omit the
$ (k+1) $-st coordinate. Let us define $ S_k $ on the algebra generators.
Here and in the sequel, it will be useful to make the following notation.
Let $ \varphi _k:\{ 1,\dots ,n\} \to \{ 1,\dots ,n+1 \} $ (where
$ 1\leq k\leq n+1 $) be the order-preserving bijection onto
$ \{ 1,\dots ,n+1 \} \setminus \{ k\} $. Set
\begin{equation}
\label{eq:51}
   S_k:=\iota _{\varphi _{k+1}} , \mbox{ on } \oH ^{\otimes n }
\end{equation}
(where $ \iota _{\varphi _{k+1}} $ is defined like  in~(\ref{eli7})), and
\begin{equation}
\label{eq:52}
   S_k(G_{ij}):=G _{\varphi _{k+1}i,\varphi _{k+1}j} \, ,
\end{equation}
for $ n\geq i>j\geq 1 $.

The face maps, $ D_0 $ and $ D_{n+1}:E_{n+1}(\oH )\to E_n(\oH ) $,
correspond to the natural inclusions,
$ L,R : F(\oX ,n)\hookrightarrow F(\oX ,n+1)$,
defined as follows. View  $ \oX $ from $ F(\oX ,n) $ as
$ X\setminus C $, where $ C $ is a $ 2m $-cube. Divide
$ C $ into two half-cubes, $ C:=C_1\cup C_2 $, pick a base-point,
$ \pt \in int(C_1) $, and view $ \oX $ from $ F(\oX ,n+1) $ as
$ X\setminus C_2 $. Then $ L(R) $ adds $ \pt $ to a given
$ n $-configuration, on the left (respectively right); see~\cite{D}.
Algebraically, we define $ D_0 $ and $ D_{n+1} $ as follows. On
$ \oH ^{\otimes n+1} $:
\begin{equation}
\label{eq:53}
   \left \{
     \begin{array}{ll}
              D_k:=\varepsilon \otimes id  , & \textrm { for }k=0 ,  \\
              D_k:=id\otimes \varepsilon  , & \textrm { for }k=n+1 .
     \end{array}
   \right .
\end{equation}
On the exterior algebra generators:
\begin{equation}
  \label{eq:54}
     \left \{
         \begin{array}{ll}
              D_0(G_{ij}):=G_{i-1,j-1}, & \textrm { for }n+1\geq i>j>1 ,\\
              D_0(G_{i1}):=0  , & \textrm { for }n+1\geq i>j=1 ,
         \end{array}
     \right .
\end{equation}
and
\begin{equation}
  \label{eq:55}
     \left \{
        \begin{array}{ll}
              D_{n+1}(G_{ij}):=G_{ij}, & \textrm { for }n+1> i>j\geq 1 ,\\
              D_{n+1}(G_{n+1,j}):=0  , & \textrm { for }n+1=i>j\geq 1 .
        \end{array}
     \right .
\end{equation}

It remains to define the face maps $ D_k:E_{n+1}(\oH )\to E_n(\oH ) $,
for $ 1\leq k\leq n $. For $ X=S^2 $, they correspond to the {\em doubling }
maps, $ \delta _k:F(\oX ,n)\rightarrow F(\oX ,n+1) $. They add to a given
$ n $-configuration, $ (x_1,\dots ,x_n) $, a new coordinate, $ x'_k $, a
copy of $ x_k $, placed to the left of $ x_k $ and close to $ x_k $; see
again~\cite{D}. Algebraically, on $\oH ^{\otimes n+1}$:

\begin{equation}
\label{eq:56}
D_k(h_1\otimes \dots \otimes h_{n+1}):=h_1\otimes \dots \otimes h_{k-1}
       \otimes (h_k\cdot h_{k+1})\otimes h_{k+2}\otimes \dots \otimes h_{n+1},
\end{equation}
and, for $ n+1\geq i\neq j\geq 1 $ :
\begin{equation}
\label{eq:57}
       \left \{
          \begin{array}{ll}
                D_k(G_{ij}):=G_{rs}, & \textrm { if } i=\varphi _kr
                                    \textrm{ and } j=\varphi _ks ;\\
                D_k(G_{kj}):=G_{ks}, & \textrm { if } j=\varphi _ks
                                    \textrm{ and } s\neq k ;\\
                D_k(G_{k,k+1}):=0  , &
          \end{array}
       \right .
\end{equation}
where $\varphi_k$ is like in~(\ref{eq:51}).

\begin{proposition}
\label{prop:drin}
The maps
\[
\{ D_k:\oH ^{\otimes n+1}[G_{ij}]\to \oH ^{\otimes n}[G_{ij}]\}_%
{0\leq k\leq n+1}
\mbox { and }
\{ S_k:\oH ^{\otimes n}[G_{ij}]\to
\oH ^{\otimes n+1}[G_{ij}]\}_{0\leq k\leq n},
\]
defined as above, induce a simplicial DBGA structure on
$ \EE (\oH ) $.
\end{proposition}

\begin{proof}
By construction, all maps are multiplicative and bihomogeneous,
on $ \oH ^{\otimes n+1}[G_{ij}] $ and $ \oH ^{\otimes n}[G_{ij}] $
respectively. Starting from definitions (\ref{eq:51})-(\ref{eq:52}),
it is straightforward to check that one has induced DBGA maps,
$ \{ S_k:E_n(\oH )\to E_{n+1}(\oH )\} _{0\leq k\leq n} $, as asserted.
Similarly, for $ D_0 $ and $ D_{n+1} $.

Let us prove this now for $ D_k $, with $ 1\leq k\leq n $. Consider first
the quotient algebra, $ \ext (G_{ij};n+1\geq i>j\geq 1) $, modulo
the Arnold relations~(\ref{one}). After uniform rescaling of degrees
(that is, putting all $ G_{ij} $'s in degree $ 1 $, instead of $ 2m-1 $),
this graded algebra becomes isomorphic to $ H^*F(\bbR ^2,n+1) $. Moreover,
the generator $ G_{ij} $ is dual to $ t^{ij} $, where $ t^{ij}\in
H_1F(\bbR ^2,n+1) $ is the homology class of the standard generator
$ \sigma ^{ij}\in \pi _1F(\bbR ^2,n+1) $ from~\cite[Proposition~3.6]{BN}.
See also~\cite{OT}. With these identifications, our formulae (\ref{eq:57}),
defining $ D_k $ on the generators $ G_{ij} $, become dual to those
from~\cite[Definition~2.9]{BN}, which express the {\em doubling } operations
on $ H_1F(\bbR ^2,n) $ in terms of the generators $ t^{ij} $. (See
\cite[\S~\S 4.4-4.5]{BN} for the importance of doubling, viewed in the
framework of {\em chord diagrams } from the theory of finite type invariants
of links.) In this way, $ D_k $ is identified with $ H^*\delta _k $ (where
$ \delta _k:F(\bbR ^2,n)\to F(\bbR ^2,n+1) $ is the geometric doubling),
on the exterior generators; in particular, $ D_k $ preserves
relations~(\ref{one}). The preservation of relations (\ref{2}) follows
immediately from definitions (\ref{eq:56})-(\ref{eq:57}), as well
as the commutation relations,
\begin{equation}
\label{eq:58}
     D_k\od G_{ij}=\od D_kG_{ij},
\end{equation}
in the first two cases from (\ref{eq:57}). Finally, $ D_k\od G_{k+1,k}=
\sum _{\alpha }(-1)^{deg(h_{\alpha })}\iota _k
(h_{\alpha }\cdot h_{\alpha }^*) $,
according to the definitions. This sum equals zero, as asserted, since
$ h_{\alpha }\cdot h_{\alpha }^*=\omega $, for all $ \alpha $, and
$ \omega =0 $ in~$ \oH $.

To finish the proof of our Proposition, it is enough to check the
simplicial identities on the algebra generators. On the exterior part,
we have identified $ D_ k $ with $ H^*\delta _k $, for $ 1\leq k\leq n $.
Formulae (\ref{eq:54})-(\ref{eq:55}) readily imply the identifications
$ D_0\equiv H^*L $ and $ D_{n+1}\equiv H^*R $ respectively. It is equally
easy to check that $ S_k\equiv H^*(pr_{k+1}) $, for $ 0\leq k\leq n $,
starting from definition (\ref{eq:52}).
Therefore, it is enough to check the {\em cosimplicial } identities,
for the continuous maps $ \{ s^k:=pr_{k+1} \} _{0\leq k\leq n} $, and
$ \{ \partial ^k:=\delta _k\} _{1\leq k\leq n}\cup \{
\partial ^0:= L, \partial ^{n+1}:=R \} $, which is routine.

Similarly, one may check the simplicial identities also
on generators coming from $ \oH $, by completely straightforward
dual calculations.
\end{proof}

\begin{remarks}
\label{rks:pess1}
{\rm
Formulae (\ref{eq:51})-(\ref{eq:57}) make sense also for
$ \{ E_n(H) \}_n $. While $ \{ S_k \} _k $ induce DBGA maps,
$ \{ S_k:E_n(H)\to E_{n+1}(H) \} _k $, the others don't, in general. Indeed,
for instance $ D_{n+1}dG_{n+1,n}=\iota _n\omega \neq 0 $, while
$ dD_{n+1}G_{n+1,n}=0 $. This is related to the fact that, for a {\em closed }
manifold $ X $, $ pr_{n+1}:F(X,n+1)\to F(X,n) $ does not need to
have a section.

Similarly, $ D_ndG_{n+1,n}=e(H)\cdot \iota _n \omega$,
while $ dD_nG_{n+1,n}=0 $. Thus, $ D_nd\neq dD_n $, when the characteristic
of $ \, \bbK $ does not divide the Euler characteristic $ e(H) $.
The topological counterpart of the potential
failure of type (\ref{eq:58}) relations, in the closed case, is the
fact that the doubling operations, $ \{ \delta _k:F(X,n)\to F(X,n+1)\} _k $,
have no natural definition, when $ X $ does not admit a nowhere zero
vector field.
\endex
}
\end{remarks}

\subsection{A 'coaction' map}
\label{subsec:act}

There is a natural topological action,
\begin{equation}
\label{eq:59}
F(\oX ,n_0)\times F(\bbR ^{2m},n_1)\times \dots \times F(\bbR ^{2m},n_r)
\stackrel {\alpha _{\TT}}\longrightarrow F(\oX ,n_0+n_1+\dots +n_r)\, ,
\end{equation}
for any partition $ \TT $ of $ n:=n_0+n_1+\dots +n_r $. Here,
$ \TT:=\{ T_j \} _{0\leq j\leq r} $, $ \{ 1,\dots ,n\} =T_0\coprod
T_1\coprod \dots \coprod T_r $, and $ n_j:=|T_j| $, for $ j=0,\dots ,r $.
For each $ j $, denote by $ \varphi ^j $ the (unique) order-preserving
bijection, $ \varphi ^j:\{ 1,\dots ,n_j\}\stackrel {\sim }\rightarrow T_j $.
To define the above action map, $ \alpha _{\TT } $, start by viewing
$ \oX $ from $ F(\oX ,n_0) $ as $ X\setminus C $, where $ C $ is a
$ 2m $-cube. Divide $ C $ into two cubes, $ C=C'\cup C'' $, and view
$ \oX $ from $ F(\oX ,n) $ as $ X\setminus C'' $.
Next, subdivide $ C' $ into $r$ cubes, $ C'=C'_1\cup \dots \cup C'_r $.
Now, let $ \underline x _0 $
be an $ n_0 $-configuration in $ \oX $, and let $ \underline x _j $
be given $ n_j $-configurations in $ \bbR ^{2m} $, for $ 1\leq j\leq r $.
Define $ \alpha _{\TT }(\underline x _0,\underline x _1,\dots ,
\underline x _r) $ to
be the $ n $-configuration in $ \oX $ obtained by  putting the
coordinates of $ \underline x _0 $ in $ X\setminus C $, on the positions
prescribed by $ \varphi ^0 $, and those of  $ \underline x _j $
($ 1\leq j\leq r $), suitably rescaled, in $ int(C'_j ) $, according to
$ \varphi ^j $-prescriptions.

Note that all spaces $ F(\bbC ^m,k) $ are formal, as follows from work by
S.~Yuzvinsky in~\cite{Y}. Their cohomology algebras were computed by F.R.
Cohen (see~\cite{C}); they have exterior generators $ \{ G_{ij}; k\geq
i>j\geq 1 \} $, in degree $ 2m-1 $, and defining relations (\ref{one}).

With these preliminaries, we may now define the algebraic analog of
(\ref{eq:59}), that is, a DBGA map,
\begin{equation}
\label{eq:510}
  q_{\TT }:(E_n(\oH ),\od )\to (E_{n_0}(\oH ),\od )\bigotimes (\bigotimes ^r
          _{j=1}H^*F(\bbR ^{2m},n_j),0),
\end{equation}
where the bigrading on $ \otimes _jH^*F(\bbR ^{2m},n_j) $ comes from
putting all $ G $-type generators in bidegree $ (2m-1,1) $, as usual.

Let us define $ q_{\TT } $ on algebra generators. For $ n\geq i\neq j \geq 1 $,
set
\begin{equation}
\label{eq:511}
       \left \{
           \begin{array}{ll}
                q_{\TT }(G_{ij}):=\iota ^k(G_{i'j'}), & \textrm{ if }
                   i=\varphi ^ki'\textrm{ and }j=\varphi ^kj'\, ;\\
                q_{\TT }(G_{ij}):=0, & \textrm{ otherwise },
           \end{array}
       \right .
\end{equation}
where $ \iota ^k $ denotes the canonical inclusion into the tensor product
of $ E_{n_0}(\oH ) $ (respectively $ H^*F(\bbR ^{2m},n_k)$), for
$ k=0 $ (respectively $ 1\leq k\leq r $). On $ \oH ^{\otimes n} $, which
is generated by $ \iota _i(h) $, where $ h\in \oH $ and $ 1\leq i\leq n $,
set:
\begin{equation}
\label{eq:512}
       \left \{
           \begin{array}{ll}
                q_{\TT }\iota _i(h):=\varepsilon (h), & \textrm{ if }
                   i=\varphi ^ki'\textrm{ and }k>0 ;\\
                q_{\TT }\iota _i(h):=\iota _{i'}(h), & \textrm{ if }
                   i=\varphi ^0i' ,
          \end{array}
       \right .
\end{equation}
where $\varepsilon $ stands as usual for the canonical augmentation of the
connected graded algebra $ \oH $.

\begin{proposition}
\label{prop:coact}
The bigraded algebra map,
$$  q_{\TT }:\oH ^{\otimes n}[G_{ij};n\geq i\neq j \geq 1]\to E_{n_0}
   (\oH )\bigotimes (\bigotimes ^r_{j=1}H^*F(\bbR ^{2m},n_j))\, ,  $$
defined by  (\ref{eq:511})-(\ref{eq:512}) above, induces a DBGA map,
as in (\ref{eq:510}).
\end{proposition}

\begin{proof}
The compatibility of $ q_{\TT } $ with the Arnold relations
(\ref{one}) from $ E_n(\oH ) $ follows from an easy analysis of the definition
given in~(\ref{eq:511}).
Likewise, (\ref{2})-compatibility is an easy consequence
of (\ref{eq:511})-(\ref{eq:512}). Finally, we have to check the
commutation of $ q_{\TT } $ with differentials, on all generators
$ G_{ij} $ , $ n\geq i\neq j \geq 1 $. In the second case from (\ref{eq:511}),
as well as in the first case (subcase $ k>0 $), this amounts to
verifying that $ q_{\TT }\oD _{ij}=0 $. This in turn follows from the
remark that, in all these cases, either $ i\notin T_0 $ or  $ j\notin T_0 $.
This implies, via (\ref{eq:512}), that  $ q_{\TT } $ sends $ \oD _{ij} $
to zero, since obviously $ (\varepsilon \otimes id)(\oD )=(id\otimes
\varepsilon )(\oD )=(\varepsilon \otimes \varepsilon )(\oD )=0 $;
see Definition~\ref{def:circ}. In the remaining case,
(\ref{eq:511})-(\ref{eq:512}) readily imply that
$ \od q_{\TT }(G_{ij})=q_{\TT }\od (G_{ij})=\oD _{i'j'} $.
\end{proof}

\begin{remark}
\label{rk:pess2}
{\rm
Again, the existence of the natural topological action (\ref{eq:59}) is
a peculiarity of the punctured case. This phenomenon is also reflected
by algebra. For example, take $ \TT=\{ T_0,T_1\} $, with
$ T_0=\{ 1,\dots ,n\} $ and $ T_1=\{ n+1\} $. Compare
(\ref{eq:53}) and (\ref{eq:55}) with (\ref{eq:511}) and (\ref{eq:512})
to infer that, for this choice of $ \TT $, $ q_{\TT }=D_{n+1} $ from
\S~\ref{subsec:simpl}. As noted in Remarks ~\ref{rks:pess1}, the definition
of  $ q_{\TT } $ would make sense also in the closed manifold case, but in
that case $ q_{\TT } $ would {\em not} commute with $ d $.
\endex
}
\end{remark}

\subsection{Connected sum}
\label{subsec:sum}

We first recall the algebraic analog of the {\em connected sum }
operation. Let $ (H,\omega _H) $ and $ (K,\omega _K) $ be two oriented
$ \bbK $-Poincar\' e duality algebras, of the same formal dimension, $ 2m $.
Define the underlying graded vector space of their connected sum  by
$ H\harp K:=\bbK\cdot 1\oplus \oH ^+\oplus \oK ^+\oplus \bbK\cdot \omega $,
with $ \omega $ in degree $ 2m $. Extend the multiplications of $ H $ and
$ K $ by setting $ h\cdot k=0 $, for $ h\in \oH ^+ $ and $ k\in \oK ^+ $.
In this way, $  (H\harp K, \omega ) $ becomes a Poincar\' e duality
algebra, endowed with two (multiplicative) canonical projections,
$ \pi _H:H\harp K \to \oH $, and $ \pi _K:H\harp K \to \oK $.

The connected sum operation for oriented manifolds provides a natural map,
$$ F(\oX,r)\times F(\oY,s)\to F(X\harp Y,r+s) , $$
which sends $ ((x_1,\dots ,x_r),(y_1,\dots ,y_s)) $ to
$ (x_1,\dots ,x_r,y_1,\dots ,y_s) $.
Let us define its DBGA analog,
\begin{equation}
\label{eq:515}
  \chi :(E_{r+s}(H\harp K),d)\to (E_r(\oH ),\od )\bigotimes (E_s(\oK ),\od ).
\end{equation}
On exterior degree zero algebra generators, set
\begin{equation}
\label{eq:513}
  \chi \mid (H\harp K)^{\otimes (r+s)}:=\pi _H^{\otimes r} \otimes
              \pi _K^{\otimes s}.
\end{equation}
On $ G $-type generators, we will define $ \chi $ by
\begin{equation}
\label{eq:514}
       \left \{
           \begin{array}{ll}
             \chi (G_{ij}):=G^H_{ij}, & \textrm{ for } r\geq i>j\geq 1;\\
             \chi (G_{r+i,r+j}):=G^K_{ij}, & \textrm{ for } s\geq i>j\geq 1;\\
             \chi (G_{ij}):=0, & \textrm{ otherwise } .
           \end{array}
       \right .
\end{equation}

\begin{proposition}
\label{prop:conn}
Formulae (\ref{eq:513})-(\ref{eq:514}) above define an induced
DBGA map, $ \chi $, as in (\ref{eq:515}).
\end{proposition}

\begin{proof}
By construction, $ \chi :(H\harp K)^{\otimes (r+s)}[G_{ij}]
\to \oH ^{\otimes r}[G^H_{ij}] \otimes \oK ^{\otimes s}[G^K_{ij}] $ is
multiplicative and bihomogeneous. The preservation of relations
(\ref{one})-(\ref{2}) is immediate.

To check the commutation relations
$$ \od \chi (G_{ij})=\chi d(G_{ij}), $$
just note that $ \Delta =\omega \otimes 1+1\otimes \omega + \oD ^H
  +\oD ^K $, by construction, and then use the definitions.
\end{proof}

{\small
\def\cprime{$'$}

\catcode`\@=11

\noindent
B.~B.: Institute of Mathematics "Simion Stoilow," P.O. Box 1--764,
RO--70700 Bucharest, Romania, \hfill\break\noindent
\hphantom{B.~B.: } email: {\tt Barbu.Berceanu@imar.ro}

\noindent
M.~M.: Mathematical Institute of the Academy, \v Zitn\'a 25, 115 67
Praha 1, The Czech Republic,\hfill\break\noindent
\hphantom{M.~M.: } email: {\tt markl@math.cas.cz}

\noindent
{\c S}.~P.: Institute of Mathematics "Simion Stoilow", P.O. Box 1--764,
RO--70700 Bucharest, Romania, \hfill\break\noindent
\hphantom{{\c S}.~P.: } email: {\tt Stefan.Papadima@imar.ro}
}

\vfill
\end{document}